%5-10-97
%9-12-96
\input amstex
\magnification=1200
\loadmsam
\loadmsbm
\loadeufm
\loadeusm
\UseAMSsymbols

\hsize=7.0truein
\hoffset=-0.35truein
\vsize=9truein
\voffset=-0.2truein

\def\leftitem#1{\item{\hbox to\parindent{\enspace#1\hfill}}}

\def\boxit#1#2{\hbox{\vrule
	\vtop{%
	\vbox{\hrule\kern#1%
	\hbox{\kern#1#2\kern#1}}%
	\kern#1\hrule}%
	\vrule}}

\def\leaderfill{\leaders\hbox to 1em{\hss.\hss}\hfill}

\parskip=\medskipamount
\document
%nopagenumbers

\input epsf

\input epsf
%centerline{\epsfbox{pic2-1.eps}}

%input psfig

%centerline{\bf }
%bigskip

\centerline{\bf Geodesic Length Functions and Teichm\"uller Spaces }
\centerline{\bf Feng Luo}

\centerline{Department of Mathematics, Rutgers University, New Brunswick, NJ 08903}

\centerline{e-mail: fluo\@math.rutgers.edu}

%footnote{Supported in part by the NSF}

{\bf Abstract}      Given a compact orientable surface with finitely
many  punctures $\Sigma$, let
$\Cal S(\Sigma)$ be the set of isotopy classes of essential  unoriented
simple closed curves in $\Sigma$. We determine a complete set of 
relations for a function
from $\Cal S(\Sigma)$ to $\bold R$ to be the geodesic length function of 
a hyperbolic metric with geodesic boundary and cusp ends on $\Sigma$. As a consequence, the Teichm\"uller
space of hyperbolic metrics with geodesic boundary and cusp ends on $\Sigma$
is reconstructed from an intrinsic $(\bold QP^1, PSL(2, \bold Z))$
structure on $\Cal S(\Sigma)$.

{\bf \S0. Introduction}

Let $\Sigma = \Sigma_{g,r}^s$ be a compact oriented  surface of genus $g$ with
$r$ boundary components and $s$ punctures, i.e., a surface of 
signature $(g,r,s)$.
The Teichm\"uller space of isotopy classes of hyperbolic metrics with geodesic
boundary and cusp ends on $\Sigma$ is denoted by $T_{g,r}^s =
T(\Sigma)$ and the isotopy
classes of essential simple closed unoriented curves in $\Sigma$ is denoted by
$\Cal S = \Cal S(\Sigma)$. A simple loop in $\Sigma$ is called 
 \it parabolic \rm if it is
homotopic into an end of $\Sigma$. The set of isotopy classes of essential
parabolic simple loops in $\Sigma$ is denoted by $P(\Sigma)$. For each $m \in
T(\Sigma)$ and $\alpha \in \Cal S(\Sigma)$, let $l_m(\alpha)$ be the length
of the geodesic representing $\alpha$ if $\alpha \notin P(\Sigma)$ and let
$l_m(\alpha) = 0$ if $\alpha \in P(\Sigma)$. The goal of the paper is to
characterize the geodesic length function $l_m$ in terms of an intrinsic
 $(\bold QP^1, PSL(2, \bold Z))$ structure on $\Cal S(\Sigma)$. 

\midspace{0.1cm}
\centerline{\epsfbox{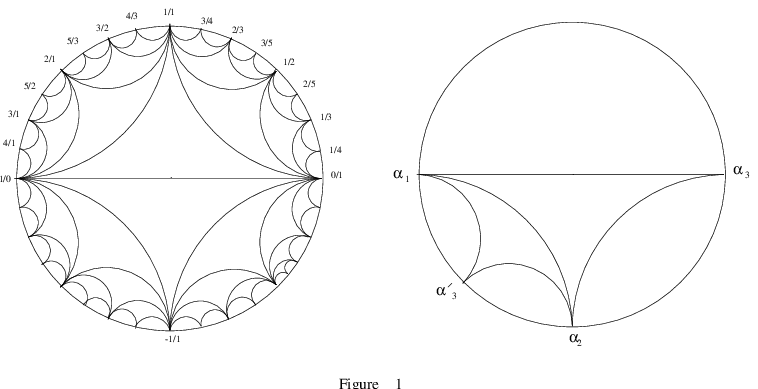}}
\midspace{0.1cm}

{\bf Theorem 1.} \it For surface $\Sigma_{g,r}^s$ of
negative Euler number, a  function $f : \Cal S(\Sigma_{g,r}^s) \to
\bold R$ is a geodesic length function if and only if $f|_{\Cal S(\Sigma')}$
is a geodesic length function for each incompressible subsurface
$\Sigma' \cong \Sigma_{1,1}^0$, $\Sigma_{0,r}^s$ $(r+s=4)$ in $\Sigma_{g,r}^s$.
Furthermore, geodesic length functions on $\Cal S(\Sigma_{1,r}^s)$ $(r+s=1)$
and $\Cal S(\Sigma_{0,r}^s)$ $ (r+s=4)$ are characterized by two polynomial
equations  $($in $\cosh( f/2)) $ in the $(\bold QP^1, PSL(2,\bold Z))$
structure on $\Cal S$. \rm

Recall that a  subsurface $\Sigma' \subset$ $\Sigma$ is \it incompressible \rm
if each essential loop in $\Sigma'$ is still essential in  $\Sigma$. 
Given two isotopy classes $\alpha$ and $\beta$  in $\Cal S(\Sigma)$, the
\it geometric intersection number \rm between $\alpha, \beta$, denoted by
I($\alpha, \beta$) is min\{$| a \cap b| | $ $a \in \alpha$ and $b \in \beta$\}
where $|a \cap b|$  is the number of points in $a \cap b$.

Given a surface $\Sigma$, let $\Cal S'(\Sigma)$ be
the set of isotopy classes of essential, non-boundary parallel non
parabolic simple loops in $\Sigma$. For surfaces $\Sigma =$
$\Sigma_{1,r}^s$ ($r+s=1$) and $\Sigma_{0,r}^s$ ($r+s =4$), 
it is well known that there exists
a bijection $\pi: \Cal S'(\Sigma) \to \bold QP^1 (= \hat \bold Q) $
so that $p'q - pq' =
\pm 1$ if and only if $I(\pi^{-1}(p/q), \pi^{-1}(p'/q')) =1$
(for $\Sigma_{1,r}^s$) and 2 (for $\Sigma_{0,r}^s$). See figure 1.
We say that three classes $\alpha_1, \alpha_2, \alpha_3$ in
$\Cal S'(\Sigma)$ form an \it ideal triangle \rm if they correspond 
to the vertices of an ideal triangle in the modular relation under the 
map $\pi$.

For the rest of the paper, we introduce the \it trace function \rm
$t_m (\alpha) = $ 2cosh$l_m(\alpha)/2$ from $\Cal S(\Sigma)$ to
 $\bold R_{\geq 2 }$.  We will  deal with the trace function 
$t_m$ instead of $l_m$.

{\bf Theorem 2.} \it (a) For surface $\Sigma_{1,r}^s$, $r+s=1$ with $b$
as the isotopy class of  the  boundary loop or the parabolic loop, 
a function $t: \Cal S \to \bold R_{\geq 2}$  is a trace
function  if and only if
the following hold.
$$ \prod_{i=1}^3 t(\alpha_i) = \sum_{i=1}^3 t^2(\alpha_i) + t(b) -2 \quad
\text{and} \tag 1$$
$$ t(\alpha_3) t(\alpha_3') =  \sum_{i=1}^2 t^2(\alpha_i) + t(b) -2$$
where $(\alpha_1, \alpha_2, \alpha_3)$ and $(\alpha_1, \alpha_2, \alpha_3')$
are distinct ideal triangles in $\Cal S'$.

(b) For surface $\Sigma_{0,r}^s$, $r+s =4$, let 
$b_1, b_2,  b_3, b_4$ be four isotopy classes of simple loops represented
by the boundary components and the parabolic loops, a function 
$t: \Cal S \to \bold R_{\geq 2}$ is a trace 
function if and only if for each ideal triangle $(\alpha_1, 
\alpha_2, \alpha_3)$ so that $(\alpha_i, b_j, b_k)$ bounds 
a $\Sigma_{0,3}^0$ in $\Sigma_{0,r}^s$ the following
hold.
$$\prod_{i=1}^3 t(\alpha_i)  = \sum_{i=1}^3 t^2(\alpha_i) + \sum_{j=1}^4
t^2(b_j) + \prod_{j=1}^4 t(b_j) + \frac{1}{2} \sum_{i=1}^3 \sum_{j=1}^4 f(\alpha_i)f(b_j)f(b_k) -4 \quad \text{and} \tag 2$$
$$ f(\alpha_3)f(\alpha_3') = 
\sum_{i=1}^2 t^2(\alpha_i) + \sum_{j=1}^4
t^2(b_j) + \prod_{j=1}^4 t(b_j) + \frac{1}{2} \sum_{i=1}^2 \sum_{j=1}^4 f(\alpha
_i)f(b_j)f(b_k) -4 $$
where ($\alpha_1$, $\alpha_2$, $\alpha_3'$) and 
 $(\alpha_1, \alpha_2, \alpha_3)$ are two distinct ideal triangles in 
$\Cal S'$.
\rm

Part (a) of theorem  2 was a result of Fricke-Klein [FK] and Keen [Ke].

Thurston's  compactification of the Teichm\"uller space  $T(\Sigma)$ (see [Bo1],
[FLP], [Th]) uses  the embedding $\tau : T(\Sigma) \to \bold R^{\Cal
S(\Sigma)}$ sending $m$ to $l_m$. Theorems 1 and 2 gives a complete
description of the image of the embedding. 

The modular relation  on $\Cal S$ is derived from an intrinsic combinatorial
structure on $\Cal S$ as. 
If two simple closed curves $a$ and
$b$ intersect at one point transversely (resp. $\alpha$, $\beta \in
 \Cal S(\Sigma)$
with I($\alpha, \beta$) = 1),
we denote it by $a \perp b$ (resp. $\alpha \perp \beta$); if two simple
closed curves $a$ and $b$ intersect at two points of different signs
transversely  and I($[a], [b]$) = 2,
we denote it by $a \perp_0 b$. In this case,  we  denote the relation
between their isotopy classes by $[a] \perp_0 [b]$.
Suppose $x$ and $y$ are two arcs in $\Sigma$ so that $x$ intersects $y$
transversely at one point. Then \it the resolution of $x \cap y$ from
$x$ to $y$ \rm is defined as follows. Take any orientation on $x$ and
use the orientation on $\Sigma$ to determine an orientation on $y$. Now
resolve the intersection point $x \cap y$ according to the orientations
(see figure 2(a)).
If $a \perp b$  or $a \perp_0 b$, we define $ab$ to be the
curve obtain  by resolving intersection points in $a \cap b$
from $a$ to $b$.  We define $\alpha \beta = [ab]$ where $a \in
 \alpha, b \in \beta$ with $|a \cap b| = I(\alpha, \beta)$.
It follows from the definition that $\alpha \beta \perp$ (resp.
$\perp_0$) $ \alpha, \beta$ if $\alpha \perp \beta$ (resp. $\alpha 
\perp_0 \beta$). Furthermore,  $\alpha(\beta \alpha) = (\alpha \beta)
\alpha = \beta$.
For surface $\Sigma =$
$\Sigma_{1,r}^s$ ($r+s=1$) and $\Sigma_{0,r}^s$ ($r+s = 4$),
three elements $\alpha_1, \alpha_2, \alpha_3$  in $\Cal S'(\Sigma)$
form an ideal triangle 
if and only if $\alpha_1 \perp \alpha_2$ and $\alpha_3 = \alpha_1 \alpha_2
$ or $\alpha_2 \alpha_1$.
In particular the two distinct ideal triangles in theorem 2 are 
$(\alpha_1, \alpha_2, \alpha_1 \alpha_2)$ and $(\alpha_1, \alpha_2, \alpha_2
\alpha_1)$.

\midspace{0.1cm}
\centerline{\epsfbox{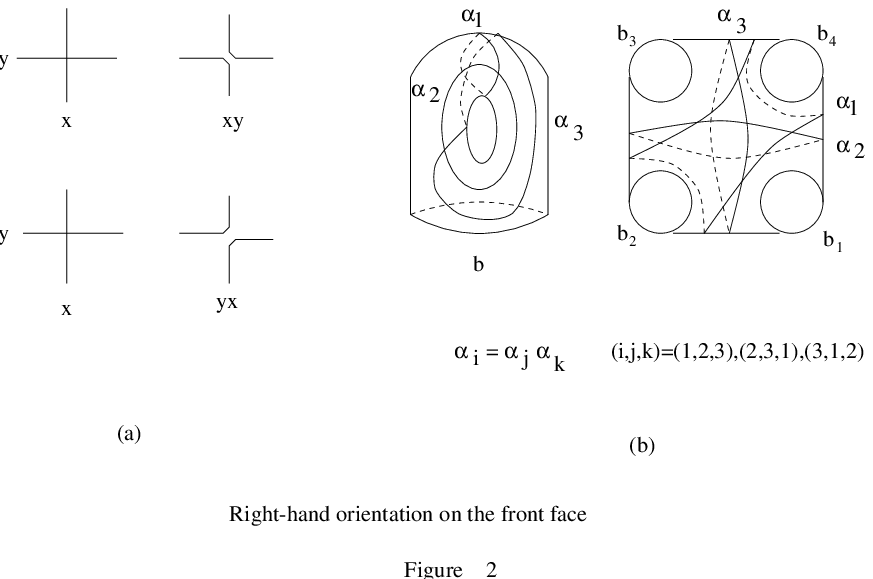}}
\midspace{0.1cm}

The relations  (1) and (2) come from trace identities for
SL(2,$\bold R$) matrices. Note that the second part of relations (1), (2)
shows that $t(\alpha_3)$ and $t(\alpha_3')$ are the two roots of
the quadratic equation (in $t(\alpha_3)$) in the first part of the relations.
Thus we obtain two more relations as follows.
$$ t(\alpha_1 \alpha_2) + t(\alpha_2 \alpha_1) = t(\alpha_1) t(\alpha_2)
\quad \text{where $\alpha_1 \perp \alpha_2$ and}$$
$$ t(\alpha_1 \alpha_2) + t(\alpha_2 \alpha_1) =
t(\alpha_1)t(\alpha_2) - t(b_i)t(b_j) -t(b_k) t(b_l) \quad \{i,j,k,l\} =
\{1,2,3,4\}$$  
where $\alpha_1 \perp_0 \alpha_2$ and $(\alpha_1 \alpha_2, b_i, b_j)$
bounds a $\Sigma_{0,3}^0.$
  
The main part of the proof of theorems is to show that
these relations are sufficient. To prove this, we  use induction on
 $|\Sigma_{g,r}^s| = 3g+r+s$.  There are two key ingredients
involved in the proof: a gluing lemma and an iteration process.

For simplicity, we describe the gluing lemma for a compact surface
$\Sigma$.  Decompose $\Sigma = X \cup Y$ where $X$ and $Y$ are
compact incompressible subsurfaces so that 
$ X \cap Y \cong \Sigma_{0,3}^0$ (see figure 3 (b), (c)).
Let the three boundary components of $X \cap Y$
be $a_1$, $a_2$ and $a_3$.
Then the gluing lemma states that for
each hyperbolic metric $m_X$ and $m_Y$ on $X$ and $Y$ respectively
so that $a_i$ are geodesics in both metrics with $l_{m_X}( a_i) = l_{m_Y}(a_i)$
(i=1,2,3), there is a hyperbolic
metric $m$ in $\Sigma$ unique up to isotopy so that the restriction
of $m$ to $X$ is isotopic to $m_X$ and the restriction of $m$ to $Y$ is
isotopic to $m_Y$.

The iteration process is derived as follows.
Given a function  $t$ on $\Cal S(\Sigma)$ satisfying the relations 
(1) and (2), using the gluing lemma and the induction hypothesis,
one constructs a hyperbolic metric on the surface
so that  $t$ and the trace of the metric coincide 
on $\Cal S(X) \cup \Cal S(Y)$.
To show that these two functions are the same on all simple 
closed curves, we observe that the second part of the relations (1) and
(2) indicates that the value of $t$ at $\beta \alpha$ is determined by the
values of $t$ on $\alpha$, $\beta$, $\alpha \beta$ and $b_i's$. 
By iterated use of the relations together with the multiplicative
structure on $\Cal S$, we show that these two functions are the same.

By the work of Thurston, the degenerations of hyperbolic
metrics become measured laminations and the corresponding projective
limits of geodesic length functions become geometric intersection
numbers. Thus, relations (1) and (2) degenerate  to universal
relations for the geometric intersection numbers. It is shown in [Lu2]
that these degenerated equations determine Thurston's measured lamination
spaces and Thurston's compactification of the Teichm\"uller spaces.

\midspace{0.1cm}
\centerline{\epsfbox{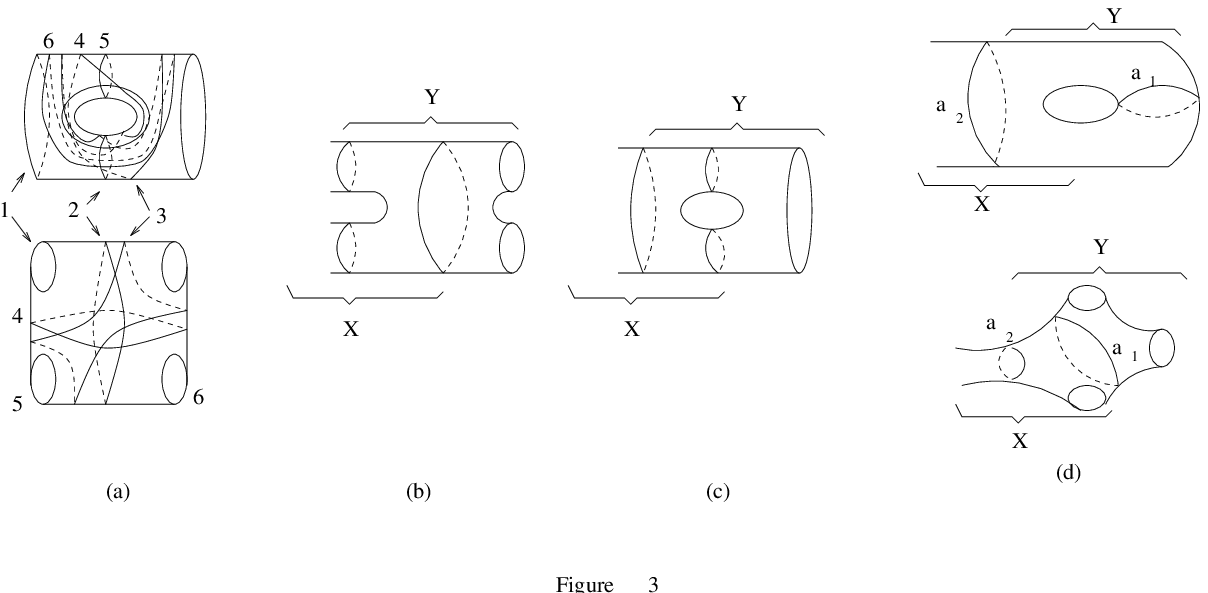}}
\midspace{0.1cm}

As another consequence of theorem 1, we consider finite dimensional
embeddings of the Teichm\"uller spaces.
Given a  subset $F$ of $\Cal S(\Sigma_{g,r}^s)$, 
let $\pi_F : T(\Sigma_{g,r}^s) \to \bold R^F$ be the map 
$\pi_F (m) = t_m |_F$.  
It is well known from the work of Fricke-Klein 
[FK] that there exists a finite set $F$ so that $\pi_F$ is an embedding.
The work of Okumura [Ok1], Schmutz [Sc],  
Sepp\"al\"a-Sorvali [SS], Sorvali [So]  show that there exists a set
$F$  consisting of $N$ (N= $6g+3r+2s -6$ if $r>0$ and  N=$6g+2s-5$ if $r=0$)
elements so that $\pi_F$ is an embedding. This number $N$ is necessary
the minimal number by a result of Wolpert  [Wo] in case $r=0$.
We shall indicate a proof of the existence of such set $F$ for compact
surface with boundary below. By theorem 2 and the gluing lemma,
it is easy to show that hyperbolic
metrics on $\Sigma_{0,4}^0$ and $\Sigma_{1,2}^0$ are determined by the
geodesic lengths of six curves as shown in figure 3(a).
Now each compact oriented surface with boundary and Euler number smaller 
than $-2$ is obtained from $\Sigma_{0,4}^0$ and $\Sigma_{1,2}^0$ by repeated 
use of gluing along 3-holed spheres (see figure 3(b), (c)). Furthermore, one
of the  subsurface used in the gluing (surface Y) is either $\Sigma_{0,4}^0$
or $\Sigma_{1,2}^0$. Thus, each time the Euler number of the resulting
surface changes by $-1$ and the number of curves needed to determine the 
hyperbolic  metric  increases by  3 (the curves $3,4,6$ in figure 3(a)
are the needed  ones and the curves $1,2,5$ are in the subsurface $X$).

The corollary below strength their result to conclude
that the image of the embedding is an explicit semi-analytic set. 
Okumura [Ok2] has also obtained the result for $s=r=0$ using a different
method.  The semi-analytic  property  in the corollary 
also follows from the work of Brumfiel [Br],  Morgan-Shalen [MS], and
Helling [He]. 

{\bf Corollary.} \it (a) 
For  surface $\Sigma_{g,r}^s$ of negative Euler number and $r >0$, there
exists a finite subset $F$ in $\Cal S(\Sigma_{g,r}^s)$ consisting of
$6g + 3r +2s -6$ elements so that the map $\pi_F : T(\Sigma_{g,r}^s)
\to \bold R^F$ is a real analytic embedding onto an open subset which is
defined by a finite set of explicit real analytic inequalities in the
coordinates of $\pi_F$.

(b) For surface $\Sigma_{g, 0}^s$ of negative Euler number, there exists
a finite subset $F$ of $\Cal S(\Sigma_{g,0}^s)$ consisting of $6g+2s -5$
elements so that $\pi_F : T(\Sigma_{g,0}^s) \to \bold R^F$ is an 
embedding whose image in $\bold R^F$ is defined
by one real analytic equation and finitely many explicit 
real analytic inequalities in the coordinates of $\pi_F$.

\rm

\midspace{0.1cm}
\centerline{\epsfbox{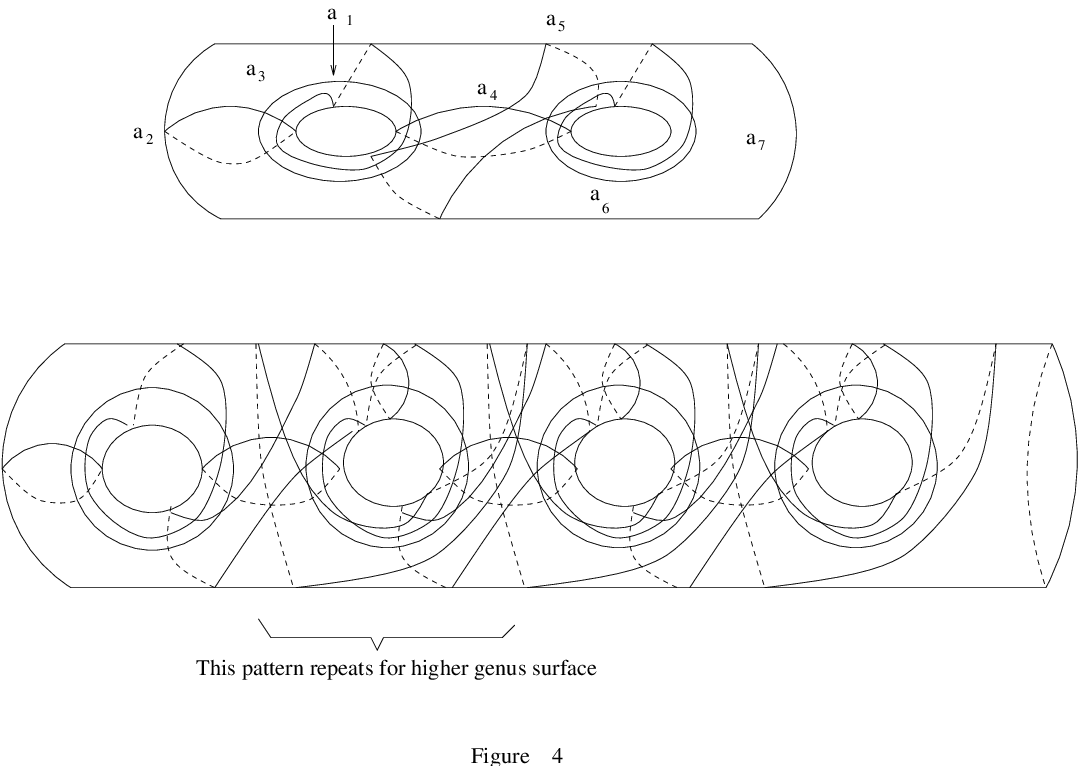}}
\midspace{0.1cm}

The inequalities and the equation in the corollary 
are given  by functions which are obtained from the coordinates of
$\pi_F$ by a finite number of algebraic operations (summation, multiplication,
and division over the rationals) and the square root operation.

Some examples of the collection $F$ and the
images of the Teichm\"uller spaces are  as follows.
For $\Sigma_{2,0}^0$, 
take $F =\{[a_1], [a_2], [a_3],[a_4], [a_5], [a_6], [a_7]\}$ as in figure 4. 
Then the  map $\pi_F$ is an embedding with image 
$\pi_F(T_{2,0}^0)$  $=\{(t_1, t_2, t_3, t_4, t_5, t_6, t_7)
\in \bold R_{>2}^7$ $|$ $t_8 > 2$, $t_9 > 2$, $t_8 = t_6t_7t_9 -t_6^2
-t_7^2 - t_9^2 + 2$, where $t_8 = t_1t_2t_3 -t_1^2 - t_2^2 -t_3^2 +2$,
and $(2 + t_2^2 + t_8) t_9^2 + 2t_2(t_4 + t_5) t_9$$
+ 2t_2^2 + t_4^2 + t_5^2 + t_8^2 + t_2^2t_8$ $-t_4t_5t_8 -4$ $= 0$\}.

The explicit equations and inequalities in the 
corollary for  the surface $\Sigma_{1,r}^s$ ($r+s=1$) are as follows.
For $\Sigma_{1,1}^0$ (resp. $\Sigma_{1,0}^1$), Keen [Ke] 
proved  that one takes
$F$ $ = \{[a_1], [a_2], [a_3]\}$ to be an ideal triangle and the
image $\pi_F (T(\Sigma_{1,1}^0))$ 
is \{$(t_1, t_2, t_3) \in \bold R_{>2}^3 |$ formula (3) holds\}.
$$ t_1 t_2 t_3 > t_1^2 + t_2^2 + t_3^3.  \tag 3$$
($\pi_F (T(\Sigma_{1,0}^1))$ =
\{$(t_1, t_2, t_3) \in \bold R_{>2}^3 |$ $t_1t_2t_3 
=$$t_1^2 + t_2^2 + t_3^2\}$).
For $\Sigma_{0,r}^s$ with
$r + s = 4$, we take the collection $F$ to be the isotopy classes 
of six curves
$b_1$, $b_2$, $b_3$, $a_{12}$, $a_{23}$, and $a_{31}$  where $[a_{ij}]$ forms
an ideal triangle and $(a_{ij}, b_i, b_j)$ bounds a $\Sigma_{0,3}$.
Then $\pi_F$ is an embedding whose image $\pi_F(T(\Sigma_{0,4}^0))$ 
is given by
\{$(t_1, t_2, t_3, t_{12}, t_{23}, t_{31}) \in \bold R_{>2}^6 |$ so that
formula (4) holds\}.
$$ t_{12}t_{23}t_{31} >  t_{12}^2 + t_{23}^2 + t_{31}^2 + t_1^2 + t_2^2
+ t_3^2 + t_{12} t_1 t_2 + t_{23} t_2 t_3 + t_{31}t_3 t_1 +  2t_1t_{23}
+ 2t_2t_{31} + 2t_{3}t_{12} + 2t_1t_2t_3.  \tag 4 $$

The organization of the paper is as follows. In \S1, we  prove
a gluing lemma and recall basic facts on discrete subgroups of
$SL(2, \bold R)$ and the spin structures on surfaces. We prove
theorem 2 in \S2. In \S3, we establish a proposition on the
multiplicative structure on $\Cal S$.  Theorem 1 is proved in \S4.
In \S5, we discuss applications.
In the main body of the paper (\S2, \S3, and \S4) we shall treat 
hyperbolic metrics without cusp ends in order to reduce the
length of the paper. No new ideas are needed for metrics with cusps. The
proofs of the theorems  1 and 2 for metrics with cups ends will 
be discussed briefly in \S5.3.

\it Acknowledgment. \rm
I would like to thank F. Bonahon  for calling my attention
to several literature. This work is supported in part by the NSF.

{\bf \S1. Preliminaries on Discrete Subgroups of SL(2,$\bold R$) }

We prove a gluing lemma in \S1.1.  Basic facts about discrete
representations of surface groups into SL(2,$\bold R$)
and spin structures
on surfaces will be recalled in \S1.2 - 1.4.

We shall use the following notations throughout the paper. Let
$\Sigma_{g,r}$ $= \Sigma_{g,r}^0$; 
$\Sigma_g$ $= \Sigma_{g,0}^0$; and $T_{g,r} = $$ T_{g,r}^0$.
We use $cl(X)$ and $int(X)$ to denote the closure and
the interior of a submanifold $X$. The isotopy class of a 
simple loop $a$ is denoted by $[a]$ and the isotopy class
of a hyperbolic metric $d$ is denoted by $[d]$.
If $f: \Cal S \to \bold R$ is a function and $a$ is a simple loop, we define
$f(a)$ to be $f([a])$. In particular, $I(a, b) = I([a], [b]) = I(a, [b])$.
A regular  neighborhood of a submanifold X is denoted by $N(X)$. Regular
neighborhoods  are always assume to be small.
All intersections of curves are assumed to be transverse.

An RC-function (compass and ruler constructible function) 
in variables $x_1$, ..., $x_n$ is a function obtained
from $1, x_1, ..., x_n$ by a finite
number of  algebraic operations  and the square root operation. The set of
RC-functions is closed under algebraic operations and compositions. Note
that $|x| = \sqrt{x^2}$ is an RC-function. An RC-function is continuous
in its natural domain and is analytic away from its singular set.

1.1.  \it A gluing lemma  \rm

First some definitions and conventions. A surface $\Sigma$ is oriented
and connected which is either $\Sigma_{g,r}^s$ or obtained from $\bar \Sigma$ =
$\Sigma_{g,r}^s$
by removing some boundary components. 
Each boundary component of $\bar \Sigma$ is called a boundary component of
$\Sigma$. A hyperbolic metric with geodesic boundary and cusp ends on $\Sigma$
is a hyperbolic metric whose completion is a hyperbolic metric  on $\bar \Sigma$
with geodesic boundary and cusp ends. Two hyperbolic metrics are \it isotopic
 \rm if there is an isometry between them which is isotopic to the identity.
The Teichm\"uller
space of hyperbolic metrics with geodesic boundary and cusp ends  on
$\Sigma$ is denoted
by $T(\Sigma)$. It is canonically isomorphic to $T(\bar \Sigma)$.

A subsurface $X$ of $\Sigma$ is \it incompressible \rm
if the inclusion map induces
a monomorphism in fundamental groups. If the subsurface is compact, then it
is incompressible if and only if each boundary component of $X$ is essential 
in $\Sigma$. A \it good incompressible \rm subsurface is an incompressible 
subsurface whose interior is a component of the complement of
a finite union of
disjoint, pairwise non-parallel, non-boundary parallel, non-parabolic simple
closed curves in $\Sigma$.  For instance, if $s$ is a non-separating
simple closed curve in $\Sigma$, then $\Sigma - s$ is a good
incompressible subsurface but $\Sigma - N(s)$ is not. 
If $X$ is an incompressible subsurface of
negative Euler number, then $int(X)$ is isotopic to a good incompressible
subsurface. For a good incompressible subsurface $X$ of $\Sigma$,
we define the \it restriction map \rm $R_{X} = R_{X}^{\Sigma} : T(\Sigma) \to T(X)$
as follows. Given $[d] \in T(\Sigma)$, there is a homeomorphism $h$ of $\Sigma$
isotopic to the identity so that the frontier of $X$,  $cl(X) -int(X)$, is
a union of geodesics in the pull back metric $h^*(d)$. We define $R_X([d])$
to be  $[ h^*(d)|_X]$. It follows from elementary hyperbolic geometry and
topology of surfaces that $R_X$ is well defined (see [CB], or [Bu]). 
Furthermore, it follows
from the definition that if $X$ is good incompressible in Y and Y is
good incompressible in Z, then $R_X^Z = R^Y_X R_Y^Z$.
The restriction map is in general not onto. For instance, if we take $X$
to be the complement of a non-separating simple closed curve in
a surface $\Sigma$ with negative Euler number, then $R_X$ is not onto.

{\bf Lemma 1.} (Gluing along a 3-holed sphere) \it  Suppose $X$ and $Y$ are
two good incompressible subsurfaces of $\Sigma$ whose union is $\Sigma$
 so that either (1) $X \cap Y$
$\cong \Sigma_{0,3}$, or (2) Y $\cong \Sigma_{1,1}$ and $X \cap Y
\cong \Sigma_{1,1} - s$ where $s$ is a non-separating simple closed curve
in int$(Y)$ (see figure 3(b), (c), (d)), or (3) $X \cap Y \cong \Sigma_{0,2}^1$ with the
punctured end in $\Sigma_{0,2}^1$ being a punctured end of $\Sigma$. Then for any two elements $m_X \in T(X)$ and $m_Y \in T(Y)$
with $R_{X \cap Y}(m_X)$ = $R_{X \cap Y}(m_Y)$, there exists a unique
element $m \in T(\Sigma)$ so that $R_X(m) = m_X$ and $R_Y(m) = m_Y$.
\rm

{\bf Proof.} To show the existence, let $d_X \in m_X$ (resp. $d_Y \in m_Y$)
be a representative so that $d_X|_{X \cap Y}$ (resp. $d_Y |_{X \cap Y}$)
has geodesic boundary and cusp ends, i.e., $R_{X \cap Y}([d_X]) = [d_X|_{X
\cap Y}]$ (resp. $R_{X \cap Y}([d_Y]) = [d_Y|_{X \cap Y}]$). Let $h$ 
$: X \cap Y \to X \cap Y$ be an isometry from $d_X |_{X \cap Y}$ to
$d_Y |_{X \cap Y}$ which is isotopic to the identity map. By the assumption
on $X$ and $Y$, we can extend $h$ to a homeomorphism $g$ of $X$ which is
isotopic to the identity. Define a hyperbolic metric $d$ on $\Sigma$ with
geodesic boundary and cusp ends as follows:  $d|_X = g^*(X)$ and $d|_Y = Y$. 
It follows from the definition that  $R_X([d]) = [d_X]$ and $R_Y([d])
= [d_Y]$. The uniqueness follows from the fact that an analytic
automorphism of a complex structure on int$(\Sigma_{0,3}$) which preserves
each end is the identity map. 
$\square$

1.2. \it Monodromy representations and spin structures \rm

Given a hyperbolic metric $d$ with geodesic boundary and cusp ends on
$\Sigma$, its monodromy is a discrete faithful representation $\rho
: \pi_1(\Sigma) \to$ PSL(2,$ \bold R)$ unique up to PGL(2,$\bold R$)
=GL(2,$\bold R )/\{\pm I\}$ conjugation so that there is an isometric
embedding $h$ from the universal cover $\tilde \Sigma$ with the
pull back metric into the hyperbolic plane $\bold H$ satisfying $h(\gamma (x)) = \rho(\gamma)(h(x))$
for all $x \in \tilde \Sigma$ and $\gamma \in \pi_1(\Sigma)$. Isotopic
metrics have the same PGL(2,$\bold R$) conjugacy class of monodromies. If the
isometric embedding $h$ is orientation preserving (resp. reversing),
we say the monodromy $\rho$
is orientation preserving (resp. reversing).
Thus each $m \in T(\Sigma)$ gives rise to two PSL(2,$\bold R)$ conjugacy
classes of monodromy representations: one preserving the orientation and
the other reversing the  orientation. Let $R(\Sigma)$ be the set of all 
such monodromy representations with the topology induced by algebraic
convergence of representations. Then $R(\Sigma)$ has two connected 
components corresponding to the two orientations. Each of the component
is a trivial principal PSL(2,$\bold R$) bundle over $T(\Sigma)$
(see [Go1], [Har], [MSi] for details). Each
representation $\rho\in R(\Sigma)$ can be lifted to a representation
$\tilde \rho: \pi_1(\Sigma) \to $ SL(2,$\bold R$) (see [Be]) 
and there are exactly $2^N$ such lifting where $N = 2g$ if $\Sigma$ has
signature (g,0,0) and  $N = 2g+r+s-1$ if $\Sigma $ has signature $(g,r,s)$
$(r+s >0)$.
Given a lifting $\tilde \rho$ of $\rho$, all other liftings are obtained as
follows. Let  \{$\gamma_1$,..., $\gamma_N$\} be a set of generators for
$\pi_1 (\Sigma)$ and $I$ a subset of $\{1,..., N\}$. Then all other
liftings are $\tilde \rho_I$ where $\tilde \rho_I ( \gamma_i) =
\tilde \rho(\gamma_i)$  if $i \in I$ and $\tilde \rho_I(\gamma_i) 
= - \tilde \rho (\gamma_i)$ if $ i \notin I$.
Let $\tilde R(\Sigma)$ be the set of all liftings of the monodromies with
the algebraic convergent topology. The representation space $\tilde R(\Sigma)$
has $2^{N+1}$ many connected components. These components are classified into
two types according to the orientation of the representations in $R(\Sigma)$.
Each component corresponds to a spin structure on the surface.
We shall recall briefly spin structures.
Let $U \Sigma$ be the unit tangent bundle over the surface $\Sigma$ with $S^1$
as a fiber. A \it spin structure \rm on $\Sigma$ is a two-fold covering space
of $U\Sigma$ so that the $S^1$-fiber does not lift. Since two-fold covering
spaces correspond to index-two subgroups of the fundamental groups,
a spin structure is the same as an epimorphism $\eta : \pi_1(U \Sigma)
\to \bold Z_2 =\{\pm 1\}$ (as a multiplicative group) so that $\eta (S^1)$ =
$-1$. Since $\bold Z_2$ is abelian, the epimorphism $\eta$ is induced by
an epimorphism (still denoted by) $\eta$ : $H_1(U \Sigma, \bold Z_2) \to \bold Z_2$ with $\eta( S^1) = -1$. Given a smooth immersed curve $c$ in $\Sigma$, let
$\bar c$ be the unit tangent vectors of $c$ in $U \Sigma$.
We define $\eta (c)$ to be $\eta ([ \bar c])$. For
instance, if $c$ bounds a disc, then $\eta( c) = -1$ and if $c$ is null
homotopic with exactly one self intersection (a figure eight), then $\eta (c) =
1$.

Johnson in [Jo] provides an algorithm to calculate $\eta (c)$ which we
summarize as follows.

{\bf Lemma 2.} ({\bf Johnson}) \it (a) Suppose \{$a_1,..., a_n\}$ and $\{ b_1,...,
b_m\}$ are two collections of disjoint simple closed curves in $\Sigma$ so
that $\Sigma_{i=1}^n [a_i] = \Sigma_{j=1}^m [b_j]$ in $H_1(\Sigma, \bold Z_2)$.
Then $\Sigma_{i=1}^n [\bar a_i] + n[S^1] = \Sigma_{j=1}^m [\bar b_j] + m[S^1]$
in $H_1( U \Sigma, \bold Z_2)$.

(b) Given $\alpha \in H_1(\Sigma, \bold Z_2)$, represent $\alpha$ as
$\Sigma_{i=1}^n [a_i]$ in $H_1(\Sigma, \bold Z_2)$ where $\{a_1,..., a_n\}$
is a collection of disjoint simple closed curves in $\Sigma$. Then
$\eta^* (\alpha)$ = $(-1)^n \Pi_{i=1}^n \eta(a_i)$ is a $\bold Z_2$-quadratic
map from $H_1(\Sigma, \bold Z_2)$ to $\bold Z_2$, i.e, $\eta^* (\alpha + \beta)
= (-1)^{<\alpha, \beta>} \eta^*(\alpha) \eta^*( \beta)$ where $<\alpha,\beta>$
 is the $\bold Z_2$-intersection number.
\rm

As a simple consequence, if $\{a_1, a_2, a_3\}$ bounds a 3-holed sphere
in $\Sigma$, then $\eta (a_1) \eta(a_2) \eta(a_3) = -1$;  if $
b$ is  the boundary of a subsurface of signature (g,1,0), then 
$\eta (b) =-1$; and  if  $a_1 \perp a_2$, then $\eta(a_1) 
\eta(a_2) \eta(a_1 a_2) = 1$.

The relationship between a lifting  $\tilde \rho\in \tilde R(\Sigma)$
of $\rho\in R(\Sigma)$ and
a spin structure is as follows. We first identify PSL(2,$\bold R$) with
$U \bold H$ by sending an isometry $g$ to $g(v_0)$ where $v_0$ is a
specified element in $U \bold H$. Under this identification, given a hyperbolic
metric with geodesic boundary and cusp ends on $\Sigma$ whose monodromy is
$\rho$, $U \Sigma$ is canonically identified with a deformation retractor
($U$(Nielsen core)) of PSL(2,$\bold R)/ \rho(\pi_1(\Sigma))$.  Let $P:$ 
SL(2,$\bold R$) $\to$ PSL(2,$\bold R$) be the canonical projection.
It is a two-fold covering map so that the $S^1$ fiber (corresponding to
PSO(2) in PSL(2,$\bold R$)) does not lift. Then $P$ induces a two-fold
covering map from SL(2,$\bold R )/\tilde \rho(\pi_1(\Sigma))$ to
PSL(2,$\bold R)/\rho(\pi_1(\Sigma))$ so that the $S^1$ fiber does not lift.
Thus we have a spin structure $\eta$ on $\Sigma$ associated to
the lifting $\tilde \rho$ of $\rho$. A simple calculation shows that
$$ \eta ( \gamma_*) = sign( tr(\tilde \rho(\gamma)), \quad \quad \gamma \in 
\pi_1(\Sigma) \tag 5$$
where $\gamma_*$ is the geodesic representative or a multiple of a parabolic
simple closed curve  in the conjugacy class of  $\gamma$.

1.3. \it Trace identities and representations of surface groups into SL(2,$\bold R$) \rm

Given three matrices $A_1$, $A_2$, $A_3$ in SL(2,$\bold C$), 
we have the following
identities on the traces of their products (see [FK], [Go2], [Ho], or [Mag]). 
The basic trace identity is
$ trA_1 A_2 + tr A_1^{-1}A_2 = trA_1 tr A_2$. By iterated use of it, one
obtains the following relations.

$$trA_1 A_2 tr A_1^{-1}A_2 = tr^2A_1 + tr^2A_2 - tr[A_1, A_2] -2. \tag 6$$
$$ tr[A_1, A_2] + 2 =  tr^2 A_1 + tr^2 A_2 + tr^2 A_1 A_2
 - trA_1 trA_2 trA_1 A_2.  \tag 7$$
$$ trA_1A_2 A_3 + trA_1 A_3 A_2 = trA_1 trA_2A_3 + trA_2 trA_3 A_1 +
trA_3trA_1A_2 - trA_1 trA_2 trA_3. \tag 8$$
$(9) \quad trA_1A_2A_3 trA_1A_3A_2 =$
$ tr^2A_1 + tr^2A_2 +$
$ tr^2A_3 + tr^2A_1A_2$
$+ tr^2 A_2A_3$
$+ tr^2 A_3 A_1$ $+$
\newline
$trA_1A_2 trA_2A_3 trA_3A_1 $
$ - trA_1trA_2 trA_1A_2 $
$- trA_2trA_3 trA_2A_3 -$
$ trA_3trA_1 trA_3A_1 -4.  $

Combining formulas (8) and (9), we see that  $tr A_1A_2A_3$ and $tr A_1A_3A_2$
are the two roots  of the quadratic equation (10) below where P and Q 
stand for the right-hand  sides of formulas (8) and (9) respectively.
$$ x^2 - P x + Q = 0.  \tag 10$$
Using  the basic  trace relation, one obtains the following (see [Ho],[CS]).

{\bf Lemma 3.} ({\bf Fricke-Klein}) \it Suppose $F_n$ is the free group on $n$
generators $\gamma_1,..., \gamma_n$. Then for each element $w$ in $F_n$, 
there is a polynomial $P_w$ with integer coefficient in $2^n-1$ variables
$x_{i_1 ... i_k}$ with $1 \leq i_1 < ...< i_k \leq n$ so that
for any representation $\rho : F_n \to$ SL(2,$\bold R$)
$$ tr \rho (w) = P_w(x_1, x_2,..., x_{i_1  ... i_k}, ..., x_{12...n})$$
where $x_{i_1 ... i_k} = tr \rho (\gamma_{i_1} ... \gamma_{i_k})$.
Furthermore,
if $\rho_1$ and $\rho_2$ are two representations with the same character
and $\rho_1(F_n)$ is not a  solvable group, then $\rho_1$ is conjugated to
$\rho_2$ by a GL(2,$\bold R$) matrix.
\rm

In particular, if $n=2$, then the three variables are $tr \rho(\gamma_1)$, $
tr  \rho(\gamma_2)$ and $tr \rho(\gamma_1 \gamma_2)$; if $n=3$, the seven
variables are $tr \rho(\gamma_i)$, and $tr \rho(\gamma_i \gamma_j)$ 
and $tr \rho (\gamma_1 \gamma_2 \gamma_3)$ where $i,j=1,2,3 $ and $i<j$.

The  discrete faithful representations of $\pi_1(\Sigma_{0,r}^s)$
($r+s =3$) and $\pi_1(\Sigma_{1,r}^s$) ($r+s=1$) which uniformize hyperbolic
structures on  $\Sigma_{0,r}^s$ ($r+s =3$) and $\Sigma_{1,r}^s$ ($r+s=1$)
are as follows.  See [GM], [Ke] for details.

For surface  $\Sigma_{0,r}^s$, $r+s =3$, we choose a set of geometric
generators $\gamma_1$ and $\gamma_2$ in $\pi_1(\Sigma_{0,r}^s)$ so that
$\gamma_1$, $\gamma_2$ and $\gamma_3 = \gamma_1 \gamma_2$ are 
represented by simple closed
curves homotopic into the three ends of $int(\Sigma_{0,r}^s)$. 
$\Sigma_{0,2}^1$ has the puncture at the end corresponding to $\gamma_3$
and $\Sigma_{0,1}^2$ has the punctures at the ends corresponding to
$\gamma_2$ and $\gamma_3$. 

{\bf Lemma 4.} ({\bf Fricke-Klein}) \it (a) If $\rho \in \tilde
 R(\Sigma_{0,r}^s)$ with $r+s = 3$,
then $tr\rho(\gamma_1) tr\rho(\gamma_2) tr \rho(\gamma_3) <0$
and $|tr \rho(\gamma_i)|$ $\geq 2$  for $i =1,2,3$  so that the  equality 
holds if and only if the corresponding end is a cusp.

(b) Given three real numbers $t_1$, $t_2$ and $t_3$ with $t_1 t_2 t_3 < 0$ and
$|t_i| > 2$ ($i=1,2,3$), there exist two elements $\rho_1$ and $\rho_2$ in $\tilde
R(\Sigma_{0,3})$ unique up to SL(2,$\bold R$) conjugation so that
$tr \rho_i (\gamma_j) = t_j$ ($i=1,2; j=1,2,3$). These two representations
are GL(2,$\bold R$) conjugated and are related by $\rho_1(\gamma_i) = \rho_2
(\gamma_i)^{-1}$. Furthermore, if $\rho(\gamma_1)$ = 
 $ \left( \matrix \lambda & 0\\ 0&  \lambda^{-1}
 \endmatrix \right)$, $\lambda > 1$, and $\rho (\gamma_2) =$ $ \left( \matrix 
a & b \\ c  & d \endmatrix \right)$, $c =1$, then $a,b$, $d$ and $\lambda$
are  real analytic RC-functions of $t_1$, $t_2$ and $t_3$ in the
domain defined by  $t_1 t_2 t_3 < 0$ and $|t_i| > 2$ ($i=1,2,3$). 

(c) Given three numbers  $t_1$, $t_2$ and $t_3$ with $t_1 t_2 t_3 < 0$ and
$|t_1| >2, |t_2| > 2$ and $|t_3| = 2$ (resp. $|t_1| >2, |t_2| = |t_3| = 
2$), there  exist two elements $\rho_1$ and $\rho_2$ in $\tilde
R(\Sigma_{0,2}^1)$  (resp. $\tilde R(\Sigma_{0,1}^2$)) unique up to SL(2,$\bold R$) conjugation so that
$tr \rho_i (\gamma_j) = t_j$ ($i=1,2; j=1,2,3$).  These two representations
are GL(2,$\bold R$) conjugated and are related by $\rho_1(\gamma_i) = \rho_2
(\gamma_i)^{-1}$.
Furthermore, if $\rho(\gamma_1)$ =  
$ \left( \matrix \lambda & 0\\ 0&  \lambda^{-1}
 \endmatrix \right)$, $\lambda > 1$, and $\rho (\gamma_2) =$ $ \left( \matrix 
a & b \\ c & d \endmatrix \right)$, $c =1$, then $a,b$, $d$ and $\lambda$
are real analytic RC-functions of $t_1$, $t_2$ and $t_3$ in the domain
defined by  $t_1 t_2 t_3 < 0$, $|t_1| >2$ and  $|t_2| > 2$.

(d) $T(\Sigma_{0,0}^3)$ consists of one point.

\rm

Note that part (a) is a consequence of lemma 2 and
formula (5). To  find the explicit expression of $a,b, d$ and $\lambda$ in terms
of $t_i's$,  see [Har], pp305.

For surface $\Sigma_{1,r}^s$ ($r+s=1$), we take a set of geometric generators
$\{\gamma_1, \gamma_2$\} in $\pi_1(\Sigma_{1,r}^s)$  so that 
they are represented by two simple closed curves $a_1$ and $a_2$ with
$a_1 \perp a_2$. The
multiplication $\gamma_3 = $$\gamma_1 \gamma_2$ is represented 
(in the free homotopy
class) by either $a_1 a_2$ or $a_2 a_1$ depending on the orientation
of the surface. The commutator $\gamma_1 \gamma_2 \gamma_1^{-1} \gamma_2^{-1}$
is represented by the simple closed curve $\partial N(a_1 \cup a_2)$
homotopic into the end of $int(\Sigma_{1,r}^s)$.

{\bf Lemma 5.} ({\bf Fricke-Klein, Keen}) \it (a)
If $\rho \in \tilde R (\Sigma_{1,r}^s)$ with $r+s=1$, then $tr \rho(\gamma_1 \gamma_2
\gamma_1^{-1} \gamma_2^{-1}) \leq -2$ so that equality holds if and only if
$s=1$. In particular,  $tr^2 \rho(\gamma_1) + tr^2 \rho(\gamma_2)
+ tr^2 \rho(\gamma_3) - tr \rho(\gamma_1) tr \rho(\gamma_2) tr \rho(\gamma_3)
\leq 0$ so that equality holds if and only if $s=1$.

(b) Give three numbers $t_i$, $i=1,2,3$ with $|t_i| > 2$ and 
$t_1^2 + t_2^2 + t_3^2 - t_1 t_2 t_3 <0$ (resp. $t_1^2 + t_2^2 +
t_3^2 - t_1t_2t_3
=0$),  there exist two representations $\rho_1$ and $\rho_2$ in $\tilde
R(\Sigma_{1,1}^0)$  (resp. $\tilde R(\Sigma_{1,0}^1)$)
unique up to SL(2,$\bold R$) conjugation so that
$tr \rho_i (\gamma_j) = t_j$ ($i=1,2; j=1,2,3$). These two representations
are GL(2,$\bold R$) conjugated and are related by $\rho_1(\gamma_i) = \rho_2
(\gamma_i)^{-1}$. Furthermore, if $\rho(\gamma_1)$ =  
$ \left( \matrix \lambda & 0\\ 0&  \lambda^{-1}
 \endmatrix \right)$, $\lambda > 1$, and $\rho (\gamma_2) =$ $ \left( \matrix
a  & b \\ c & d \endmatrix \right)$, $c =1$, then $a,b,c$, $d$ and $\lambda$
are real analytic RC-functions of $t_1$, $t_2$, $t_3$.

\rm

The first part of the lemma also follows lemma 2 and
formula (5). 
Below is a  proof of part (b) (known to J. Gilman).
By lemma 3, it suffices to show the existence of $\rho \in \tilde
R(\Sigma_{1,r}^s)$ with  $tr(\rho(\gamma_j)) = t_j$, $j=1,2,3$. We first
construct three points $A_1$, $A_2$ and $A_3$ in
$\bold H$ so that their pairwise hyperbolic distance $d(A_i, A_j)$ 
is determined by 2cosh$d(A_i, A_j)$/4 = $|t_k|$ where $i \neq j \neq k \neq i$.
That the pairwise distances satisfy the triangular inequalities follows
from the given condition on $t_i's$. Let $h_{A_i}$ be the hyperbolic 
isometry which rotates by degree $\pi$ at the point $A_i$ (a half-turn).
Then $h_{A_i}h_{A_j}$ ($i \neq j$) is a hyperbolic isometry so that the
absolute value of its trace is $|t_k|$ by the construction ($k \neq i,j$). 
Furthermore, $tr(h_{A_1} h_{A_2} h_{A_3})^2 = tr[h_{A_1}h_{A_2}, h_{A_3}h_{A_1}]$ = $ t_1^2 + t_2^2 + t_3^2 - t_1 t_2 t_3  -2$ which is at most $-2$. Thus
the isometry $h_{A_1} h_{A_2} h_{A_3}$ has a fixed point $p$ at the circle at
the infinity of $\bold H$. By the construction, the three vertices
of the triangle $A_1$, $A_2$ and $A_3$
are in the three sides of the ideal hyperbolic triangle  $\Delta$ with
vertices $p,$ $ h_{A_3}(p)$ and $h_{A_2} h_{A_3}(p)$. The four components
of the complement of the ideal quadrilateral $\Delta \cup  h_{A_3}(\Delta)$
give rise to a Schottky condition for the group $<h_{A_1}h_{A_3}, h_{A_3}
h_{A_2} >$. Thus by Poincar\'e polyhedron theorem, the group
 $< h_{A_1}h_{A_3}, h_{A_3}
h_{A_2} >$ uniformizes either $\Sigma_{1,1}^0$ or $\Sigma_{1,0}^1$ so that
the geodesics of $h_{A_1} h_{A_3}$ and $h_{A_3} h_{A_2}$ are
simple closed curves intersecting at one point. Let $Y$ be the lifting
of $h_{A_1} h_{A_3}$ to SL(2,$\bold R$) with $t_2 trY >0$ and $X$
be the  lifting of $h_{A_3} h_{A_2}$ to SL(2,$\bold R$)  with $t_1 trX >0$.
Then $trX = t_1$ and $trY = t_2$ and $tr(XY) = t_3$ due to the spin
structure. This finishes the proof.

{\bf \S2. Proof of Theorem 2 }

Given a hyperbolic metric $m$ on $\Sigma$ and a monodromy $\rho
\in \tilde R(\Sigma)$ of the metric $m$, we have $t_m(x) = |tr(\rho(x))|$ 
where $x$ is the homotopy class of a loop.

2.1. \it Proof of theorem 2 for $\Sigma_{1,1}$ \rm

To show that condition (1) in part (a) is necessary, take  three
classes $\alpha_1, \alpha_2, \alpha_3$ forming an ideal triangle in
$\Cal S$. Choose $\gamma_1, \gamma_2 \in \pi_1(\Sigma)$ so that
the homotopy classes  $\gamma_1, \gamma_2, \gamma_1\gamma_2$  and $\gamma_1^{-1}
\gamma_2$ represent $\alpha_1, \alpha_2, \alpha_3$ and $\alpha_3'$
respectively. If $t_m$ is a trace function
corresponding to a monodromy $\rho \in \tilde R(\Sigma_{0,4})$,
then condition (1) follows from the trace identities (6), (7) and lemma
5 where $A_i = \rho(\gamma_i)$.

To show that condition (1) is also sufficient, we note that the
modular relation implies that the value of $t$ is determined by
$t$ on $\{\alpha_1, \alpha_2, \alpha_3\}$ where $\alpha_i$'s
form an ideal triangle. Now since $t(b) \geq 2$,  by condition (1),
$t_i = t(\alpha_i)$ satisfies the inequalities in lemma 5. By lemma
5, we construct a hyperbolic metric $m$ so that $t_m(\alpha_i) =t_i$. Thus,
$t = t_m$ on $\Cal S$ by the modular relation.

The proof of theorem 2 for $\Sigma_{0,4}$ is in the same spirit,
but technically is more complicated.

2.2. \it Necessity of condition (2) in theorem 2 \rm

Given three classes $\alpha_1, \alpha_2, \alpha_3$ forming an ideal 
triangle in $\Cal S$, we take $a_{ij} \in \alpha_k$, $(i,j,k)=(1,2,3),
(2,3,1),(3,1,2)$ so that $|a_{ij} \cap a_{jk}| =2$. Without loss of 
generality, we may assume that $(a_{ij}, b_i, b_j)$ bounds $\Sigma_{0,3}$ in $\Sigma_{0,4}$.
Choose in $\Sigma_{0,4}$ a set of generators \{$A_1, A_2, A_3\}$ for $\pi_1(\Sigma_{0,4})$ as in figure 5 (multiplication of loops in
$\pi_1$ starts from left to right) so that (1)
the boundary  components $b_1$, $b_2, b_3,$ and $b_4$ of $\Sigma_{0,4}$ are
homotopic to representatives in $A_1$, $A_2$, $A_3$ and $A_1A_2A_3$
respectively; (2)
the curves $a_{12}, a_{23}$, and $a_{31}$ are homotopic to representatives  in
$A_1 A_2$, $A_2 A_3$, and $A_3A_1$ respectively; and (3) 
the generators are symmetric with respect to a $\bold Z_3$
action on $\Sigma_{0,4}$ preserving $b_4$ (figure 5(e)). 

\midspace{0.1cm}
\centerline{\epsfbox{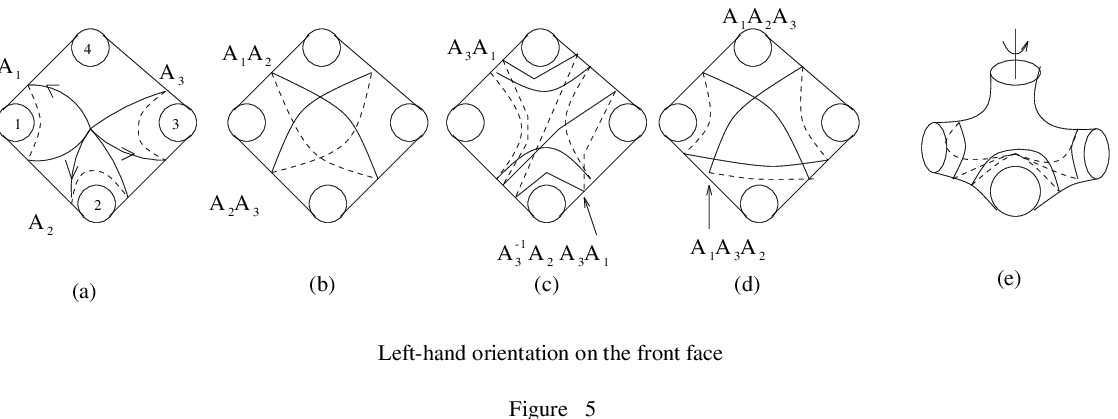}}
\midspace{0.1cm}

Given $\rho \in \tilde
R(\Sigma_{0,4})$ representing the monodromy of a hyperbolic metric
$m$, we shall identify $\rho (A)$ with $A$ for $A \in \pi_1(
\Sigma_{0,4})$ for simplicity in this section. Thus $A_i's$ are
SL(2,$\bold R)$ matrices. By choosing a different lifting if necessary, we may
assume that $tr A_i <0$ ($i=1,2,3$). By lemma 4, 
$tr A_iA_j < 0$ ($i \neq j$), and $tr A_1A_2 A_3 < 0$. Then the first
equation in condition  (2) is given by trace identity (10).
To see the second equation
(which is the statement that $f(\alpha_3), f(\alpha_3')$ are the two roots
in the first equation), we shall derive the equivalent equation 
$$ t_m(\alpha_1 \alpha_2) + t_m(\alpha_2 \alpha_1) =
t_m(\alpha_1)t_m(\alpha_2) - t(b_1)t_m(b_2) -t_m(b_3) t_m(b_4).$$
To see this, we note that $\alpha_1 \alpha_2 = \alpha_3$  and
$\alpha_2 \alpha_1$ are represented by $A_1A_2$ and $A_3^{-1}A_2A_3A_1$
respectively.  Furthermore, by lemma 4, $tr(A_3^{-1}A_2A_3 A_1) < 0$.
Thus the above formula is a consequence of
the trace identity: $tr(A_3^{-1}A_2A_3A_1) + tr(A_1A_2)
= tr(A_1)tr(A_2) + tr(A_3)tr(A_1A_2A_3) - tr(A_2A_3)tr(A_3A_1)$.

We shall write the first equation in condition (2) 
(i.e., equation (10)) explicitly as follows. Let
$t_i = t_m(b_i)$ and $t_{ij} = t_m(\alpha_k)$.
Then formulas
(8) and (9) become:
$$ -t_{4} + tr(A_1A_3A_2) = t_1t_{23} + t_2 t_{31} + t_3t_{12} + t_1t_2 t_3.
\tag 11$$
$$ -t_{4}tr(A_1A_3A_2) = t_1^2 + t_2^2 + t_3^2 + t_{12}^2 + t_{23}^2
+ t_{31}^2 + t_1t_2t_{12} + t_2t_3t_{23} + t_3t_1t_{31} -t_{12}t_{23}t_{31} -4.
\tag 12$$
Thus equation (10) becomes
$$ t_{4}^2 + t_{4} ( t_1 t_{23} + t_2 t_{31} + t_3 t_{12} +
t_1t_2t_3) + \sum_{i=1}^3 t^2_i + \sum_{(i,j) \in I}t_{i j}^2  + t_{i}t_{j}t
_{ij} -4 - t_{12}t_{23}t_{31} =0 \tag 13$$
where I = \{(1,2), (2,3), (3,1)\}.
As a quadratic equation in $-t_4 = tr(A_1A_2A_3)$,
it becomes $x^2 - Px + Q =0$ where
$P >0$ and (thus) $Q <0$. This implies that the equation has two real
roots of different signs and $-t_4$ is the negative root, i.e.,
$$ t_4 =  (-P + \sqrt{ P^2 - 4Q})/2.   \tag 14$$
In particular, the  number $t_4$ is determined by the  rest of the
six numbers.  Since $t_4 > 2$, we obtain the (equivalent) condition that
$-Q > 2P+4$ which is exactly condition (4).  Conversely, if
$-Q > 2P + 4$ and $P>0$, then $t_4 > 2$.

\it Remark 2.1. \rm We have shown that  each hyperbolic metric $m$ on
$\Sigma_{0,4}$ is determined by its lengths on six curves $\{a_{ij}, b_1, b_2,
b_3\}$.  This  was first observed by
Schmutz ([Sc], lemma 2).

2.3. \it Sufficiency of condition (2) in theorem 2\rm

We use the same notations as in \S2.2.
Given a function $t : \Cal S(\Sigma_{0,4} ) \to \bold R_{>2}$ satisfying
condition (2), we note that the
modular relation implies that the values of $t$ is determined by
$t$ on $\{\alpha_1, \alpha_2, \alpha_3, b_1, b_2, b_3, b_4\}$ 
where $\alpha_i$'s form an ideal triangle.
Thus, it suffices to find
$\rho \in \tilde R(\Sigma_{0,4})$ so that $tr \rho(A_i) = - t(b_i)$,
$tr \rho(A_iA_j) = -t(\alpha_k)$ and $tr \rho(A_1A_2A_3) =  -t(b_4)$.

Let $t_i = t(b_i)$ ($i=1,2,3,4$) and  $t_{ij}= t(\alpha_k)$, Then 
$t_i, t_{ij } \in \bold R_{>2}$ and equation (13) holds. By the 
remark in the last paragraph, this is the same as assuming condition (4)
holds for $t_1, t_2, t_3$ and $t_{ij}$.
We shall first construct three matrices $A_i$ ($i=1,2,3$) in SL(2,$\bold R)$ 
so that $tr A_i = -t_i$, $tr(A_i A_j) = -t_{ij}$ and furthermore 
$trA_1A_2 A_3 <-2$. Then we show that $tr(A_1A_2A_3) = -t_4$ and
 the corresponding representation $\rho$ is in $\tilde R(\Sigma_{0,4})$. 

Since conditions (13) and (4) is symmetric 
in $t_{12}$, $t_{23}$ and $t_{31}$ and the
set of generators $A_1$, $A_2$, and $A_3$ are also symmetric, we may assume 
without loss of generality that $t_{23} =$ max($t_{12}, t_{23}, t_{31})$.

To solve $tr A_i = -t_i$ and  $tr(A_i A_j) = -t_{ij}$, let
$A_1$ = $ \left( \matrix   x    &  y    \\  z     &   w    \endmatrix \right)$,
$A_2$ =$ \left( \matrix   a   &   b    \\   c   &   d    \endmatrix \right)$,
$A_3$ = $ \left( \matrix  -\lambda d    &  \lambda^{-1} b \\ \lambda c &  - \lambda ^{-1} a    \endmatrix \right)$ be SL(2,$\bold R$) matrices. We have,
 $$A_2A_3 =  \left( \matrix   - \lambda   &  0     \\ 0   &  -\lambda ^{-1} \endmatrix \right),$$
$$ A_1A_2 = \left( \matrix  ax+cy      &  *     \\  *    &  bz + dw    \endmatrix \right),$$
$$A_3A_1 =  \left( \matrix  -\lambda dx + \lambda^{-1}bz    &  *     \\  * & \lambda cy  - \lambda^{-1} aw   \endmatrix \right),$$
$$A_1A_2A_3 =  \left( \matrix   - \lambda x    &  *     \\  *     &  -\lambda^{-1} w \endmatrix \right).$$

By the condition $tr A_i = -t_i$ and $tr(A_i A_j) = -t_{ij}$, 
we obtain a system of quadratic and linear equations in $a,b,c,d, x,y,z,w$ and $\lambda$ as
follows.
$$
\gather
 a+ d = -t_2.  \tag E1 \\
\lambda^{-1} a + \lambda d = t_3.  \tag E2 \\
 \lambda + \lambda^{-1} = t_{23}. \tag E3 \\
 ad -bc = 1.  \tag E4 \\
 x + w = -t_1.  \tag E5 \\
 ax + cy + bz + dw = -t_{12}.  \tag E6 \\
-\lambda dx + \lambda cy + \lambda^{-1} bz - \lambda^{-1} aw = - t_{31}. \tag 
E7 \\
 xw -yz = 1.  \tag E8 \\
\endgather
$$
By (E3), $\lambda$ is a positive real number not equal to $1$ and is determined up
to reciprocal. Let us fix $\lambda > 1$. By (E1) and (E2), we have
$ a = -(\lambda t_2 + t_3)/(\lambda - \lambda^{-1})$
and
$ d = (\lambda^{-1} t_2 + t_3)/(\lambda - \lambda^{-1})$.
Thus $ad <0$ and $bc = ad -1 <0$. Fix $c=1$. We obtain a set of solutions in
$a, b, c, d$ and $\lambda$ which are real analytic RC-functions in $t_i's$ and $t_{ij}'s$.
We now claim that there are solutions for $x,y,z,$ and 
$w$ satisfying (E5)-(E8) in
the complex number field $\bold C$. Indeed, by (E6) and (E7), we express $y$ and $z$
in terms of $x$ and $w$ as follows.
$y = ( t_3x - \lambda^{-1}t_2w + \lambda^{-1} t_{12} - t_{31})/(c(\lambda - \lambda^{-1})
)$ and
$$ z = (t_2 \lambda x - t_3 w - \lambda t_{12} + t_{31})/(b(\lambda - \lambda^{-1})). \tag E9$$
Using (E5), we have $w = -x -t_1$. Thus,
$y = (\lambda^{-1} t_2 + t_3) x /(c(\lambda - \lambda^{-1})) + const$ and
$z = (\lambda t_2 + t_3)x /(b(\lambda - \lambda^{-1})) + const.$ Substitute
these
new equations and $w = -x -t_1$ into (E8). We obtain a quadratic equation in $x$
whose leading coefficient (after a simple calculation) is  $1/(bc ) \neq 0$. Thus
there is a solution for $x$ in $\bold C$. This implies the existence of
solutions for $y$, $z$ and $w$ in $\bold C$. 

We next claim that $x$, $y$, $z$, and $w$ are real numbers, i.e.,
 $A_1$ is in SL(2,$\bold R$). 
Indeed,  the quadratic
equation (in $-t_4$) (13) $x^2 -P x + Q = 0$ has two real 
roots of different signs.  By (13), 
both $trA_1A_2A_3$ and $tr A_1A_3A_2$ are solutions of the equation
Thus  $trA_1A_2A_3$  is
a real number. But $trA_1A_2A_3$ $=-\lambda x - \lambda^{-1} w$. This
together with equation (E5) shows that both $x$ and $w$ are real
numbers. Thus $y$ and $z$ are real numbers as well.

Now by choosing a different set of solution if necessary, we may assume that
$trA_1A_2A_3$ is the negative root $-t_{4}$ of the equation $t^2 -P t +
Q =0$, i.e.,
$$ \lambda x + \lambda^{-1} w = t_{4}.   \tag E10$$
Indeed, if $trA_1A_2A_3$ is the positive root, we use the new set of solution
$(A_1^{-1}, A_2^{-1}, A_3^{-1})$ to the equations $trX_i = -t_i$
and $trX_iX_j = -t_{ij}$ instead of $(A_1, A_2, A_3)$ and use the fact that
$trA_1^{-1} A_2^{-1}A_3^{-1} = trA_1 A_3 A_2$.

By the proof of above, we see that the solution $a,b,c,d, x,y,z,w$ and 
$\lambda$ are real analytic
RC-functions in $t_i$ and $t_{ij}$ (i,j=1,2,3).

By condition (4), the negative root $tr A_1A_2A_3$ is less than $-2$, i.e.,
$t_4 > 2$. Thus both
representations of $\pi_1(\Sigma_{0,3})$ (in term of the pair of
matrices) given by $<A_1^{-1}, A_1A_2A_3>$ and $<A_2, A_3>$ are in 
$\tilde R(\Sigma_{0,3})$ by lemma 4. Furthermore,  these two
group share a common generator $A_1^{-1}(A_1A_2A_3) = A_2A_3$. 
To apply the Maskit combination theorem [Ma] 
to amalgamate these two groups, we need
to verify that the Nielsen convex cores for the two groups 
$<A_1^{-1}, A_1A_2A_3>$ and $<A_2, A_3>$ in $\bold H$ lie in the 
different sides of the axis of $A_2A_3$. The
following lemma characterizes the side of the axis which contains  the
Nielsen core.

{\bf Lemma 6.} \it Suppose $X$  =$\left(\matrix  -\lambda    &  0  \\  0  &  -\lambda ^{-1} \endmatrix \right)$ and $Y$ $=$ $\left(\matrix   a  &   b  \\  c  &  d \endmatrix \right)$
are SL(2,$\bold R)$ matrices so that $trX  <-2$, $trY \leq -2$, $trXY \leq -2$. Then
the side of the axis of $X$ which contains  the Nielsen convex core
for the discrete group $<X, Y>$ is $\{(x,y)| x>0, y>0\}$
if  and only if $c (\lambda - \lambda^{-1}) > 0$.

\rm

{\bf Proof.} Let $trX = -t_1$, $trY = -t_2$, $trXY = -t_3$ with $t_1 > 2$
and $t_2, t_3  \geq 2$. Then we have $\lambda + \lambda^{-1} = t_1$, $a+d = -
t_2$,  $\lambda a + \lambda ^{-1} d = t_3$, and $ad -bc = 1$.  We solve for
$a$ and $d$ and obtain: $a = (\lambda^{-1}t_2 + t_3)/(\lambda - \lambda^{-1})$
and $d = -(\lambda t_2 + t_3)/(\lambda - \lambda^{-1})$.
The fixed points
$r_1$ and $r_2$ of $Y$ at the circle at the infinite of $\bold H$ are
the roots of the equation $c t^2 + (d-a)t  -b = 0$. In particular,
$r_1 + r_2 = -(d-a)/c$ which is $(t_1t_2 + 2t_3)/(c(
\lambda -\lambda^{-1}))$.  Since the fixed points $r_1$ and $r_2$ are
in the Nielsen core, the result follows.
$\square$

Now to finish the proof, we verify the side condition by
taking $X = A_2A_3$, and $Y= A_1$ for
the group  $<A_1^{-1}, A_1A_2A_3>$, and taking  $X= A_2A_3$,
$Y$ = $A_2^{-1}$ for $<A_2, A_3>$. Thus
it suffices to show $-zc <0$, or the same $zb <0$.

By (E5) and (E10), we have $x = (\lambda^{-1} t_1 + t_{4})/(\lambda - \lambda^{-1})$ and $w = -( \lambda t_1 + t_{4})/(\lambda - \lambda^{-1})$. Substitute
them into (E9) and simplify it, we have,
$$bz(\lambda - \lambda^{-1})^2 = t_2 \lambda(\lambda^{-1} t_1 + t_{4})
+ t_3(\lambda t_1 + t_{4}) - (\lambda - \lambda^{-1}) \lambda t_{12} +
(\lambda - \lambda^{-1}) t_{31}.$$
By (E3), we replace $\lambda^2$ by $\lambda t_{23} -1$ and $\lambda^{-1}$
by $t_{23} - \lambda$ in the above equation and obtain,
$$ bz(\lambda - \lambda^{-1})^2 = \lambda (t_1t_3 + 2t_{31} + t_2 t_{4}
-t_{23}t_{12}) + (t_1 t_2 + 2t_{12} + t_3 t_{4} - t_{23} t_{31}).$$
We claim that under the condition $t_{23} = $max$(t_{12}, t_{23}, t_{31})$
and  equation (13)
both $t_1t_3 + 2t_{31} + t_2 t_{4}
-t_{23}t_{12}$ and $t_1 t_2 + 2t_{12} + t_3 t_{4} - t_{23} t_{31}$
are negative. Indeed, since 
$t_{23}$ = max($t_{12}, t_{23}, t_{31}$), and $t_i$, $t_{ij}$ are
at least $2$, by equation (13),  we have,
$$ t_{12}t_{23}t_{31} > t_1t_3t_{31} + t_{31}^2 + t_{23}^2 + t_2 t_{31}t_{4}$$
$$ \geq  t_1t_3t_{31} + 2t_{31}^2 + t_2 t_{31}t_{4}$$
$$= t_{31}( t_1t_3 + 2t_{31} + t_2t_{4}).$$
This shows  $t_1t_3 + 2t_{31} + t_2 t_{4} -t_{23}t_{12} <0$. The other
inequality follows by the same argument since the 
inequality is obtained from the previous one by interchanging
the indices $2$ and $3$.
$\square$

The proof shows that all the entries of the matrices $A_1, A_2, A_3$
are RC functions in $t_i, t_{ij}$ where $i=1,2,3, (i,j) = (1,2),(2,3), (3,1)$.

{\bf Corollary 2.1.} \it For surface $\Sigma_{0,4}$ with $\partial
\Sigma_{0,4} = b_1 \cup b_2  \cup b_3$$\cup b_4$, let $F$
 =\{$ [a_{12}], [a_{23}], [a_{31}], b_1, b_2,$ \newline$ b_3$\} so that
$[a_{ij}]$ forms an ideal triangle and $(a_{ij}, b_i, b_j)$ bounds
a $\Sigma_{0,3}$. Then the map $\pi_F$: $T(\Sigma_{0,4})
\to \bold R^6$ is an embedding so that
its image is given by \{$(t_1, t_2, t_3, t_{12}, t_{23}, t_{31}) \in
 \bold R_{>2}^6 |$ formula (4) holds\}. Furthermore,
there exits a continuous function $f$: $T(\Sigma_{0,4})
\to \tilde R(\Sigma_{0,4})$ sending $m \in T(\Sigma_{0,4})$
a representation $f(m)$ which is a lifting of a monodromy of $m$ so that
the entries of the matrices $f(m)(\gamma)$ are real analytic
RC-functions of the
coordinates of $\pi_F(m)$, for each $\gamma \in \Cal S(\Sigma_{0,4})$.
\rm

\it Remark 2.2. \rm  The above proof works for hyperbolic metrics with cusp
ends as well since  lemmas 3, 6 and Maskit combination theorem
still hold. In particular, we obtain the following
parametrization of the Teichm\"uller space of $T_{0,0}^4$ by the
geodesic lengths $t_{12}, t_{23}$ and $t_{31}$ (other
variables $t_1$, $t_2$, $t_3$ and $t_{4}$ are $2$). 
Take $F $ =\{$ [a_{12}], [a_{23}], [a_{31}]$\}. Then the image of the
embedding $\pi_F$ of $T_{0,0}^4$ is 
\{$(t_{12}, t_{23},$$ t_{31})$$\in$ \newline
$\bold R_{>2}^3 | $ $t_{12}t_{23}t_{31} = t_{12}^2 + t_{23}^2 
$$+ t_{31}^2 + 8t_{12} + 8t_{23} + 8t_{31} + 28\}$.

{\bf \S3. A Combinatorical Structure on the Set of Isotopy Classes of 
Simple Closed Curves}

We introduce the following notation for convenience. If $\alpha
\perp_0 \beta$ (resp. $\alpha \perp \beta$),
then $\partial N(\alpha \cup \beta)$ denotes the union of the
isotopy classes of  four boundary components of $N(a \cup b)$ where
$a \in \alpha$, $b \in \beta$ with $|a \cap b| = I(a, b)$ 
(resp. $N(\alpha \cup \beta) = [\partial N(a \cup b)]$).

The goal of this section is to prove the following proposition.

{\bf Proposition 1. } \it (a) Given a set of disjoint simple closed curves and
proper arcs $\{c_1, ..., c_n\}$ in a  compact oriented surface
$\Sigma$, let $G_0 = \{ \alpha \in \Cal S(\Sigma) |$ $I(\alpha, [c_i])
\leq 2$  so that for each index $i$, if 
 equality holds  then the two points of intersection have
different signs\}. Then $\Cal S(\Sigma) = \cup_{i=0}^{\infty}G_i$ where
$G_{i+1} = G_i \cup \{ \alpha |$  $\alpha = \beta \gamma $ where either
 (1) $\beta \perp \gamma$, and $\beta$, $\gamma$,
$\gamma \beta$ are in $G_i$, or (2)
$\beta \perp_0 \gamma$, and  $\beta$, $\gamma$, $\gamma \beta$,
and each component of $\partial N(\beta \cup \gamma)$ are  in $G_i$\}.

(b) Under the same assumption as in (a), if $f$ is a function defined 
on $\Cal S(\Sigma)$ so that
(1) $f(\alpha \beta)$ is determined by $f(\alpha)$,
$f(\beta)$, and $f(\beta \alpha)$  whenever  $\alpha \perp \beta$,
and (2)  $f(\alpha \beta)$ is determined by $f(\alpha)$, $f(\beta)$,
$f(\beta \alpha)$, and $f(\gamma_i)$ ($i=1,2,3,4$) whenever
$\alpha \perp_0 \beta$ with $\partial N(\alpha \cup \beta)
= \cup _{i=1}^4 \gamma_i$,  then $f$ is determined by $f|_{G_0}$.
 
\rm

Part (b) of the proposition follows
from part (a).
The proof of  part (a) of the proposition  is a simple application of
the lemma below by induction on the number max$\{I(\alpha, [c_i])| i=1,..., n$\}
for $\alpha \in \Cal S(\Sigma)$. This lemma is motivated by lemma 2 in [Li]. 

{\bf Lemma 7.}  \it Suppose $a$  is a simple closed
curve and $b$ is either a simple closed curve or an arc so that
either $ I(a, b) = |a \cap b|  \geq 3$ 
or $a$ intersects $b$ at two points of the same intersection signs.
Let $ \{c_1, .., c_n\}$ be a collection of disjoint simple closed 
curves or arcs so that
$int(b) \cap c_i = \emptyset$ for all $i = 1,..., n$.  Then there
exist two simple closed curves $p_1$ and $p_2$ in $N(a \cup b)$ so that

(1) $a = p_1 p_2$ where either  $p_1 \perp  p_2$ or $p_1 \perp_0 p_2$,

(2) $| p_i \cap b| < |a \cap b|$, $|p_2 p_1 \cap b | < |a \cap b|$,
$|p_i \cap c_j| \leq |a \cap c_j|$ and $|p_2 p_1 \cap c_j| \leq |a \cap c_j|$
for $i=1,2$ and $j =1,2,...,n$, and,

(3) if $p_1 \perp_0 p_2$, there are four simple closed curves $d_1$, $d_2$,
$d_3$, and $d_4$ isotopic to four boundary components of $N(p_1
\cup p_2)$ so that $|d_i \cap b| < |a \cap b|$ and $|d_i \cap c_j|
\leq |a \cap c_j|$ for $i = 1,2,3,4$, and $j = 1,..., n$.

\rm

{\bf Proof}. We need to consider two cases.

\it Case 1. \rm There exist two adjacent intersection points $x$ and $y$
in $b$ which have the same intersection signs (see figure 6). 
Let $c$ be an arc in $b$ joining $x$ and $y$ so that $int(c) \cap a = 
\emptyset$. Then 
the curves $p_1$ and $p_2$ as shown in figure 6 (with the right-hand
orientation on the surface)
satisfy $p_1 \perp p_2$ and all conditions in the lemma.

\midspace{0.1cm}
\centerline{\epsfbox{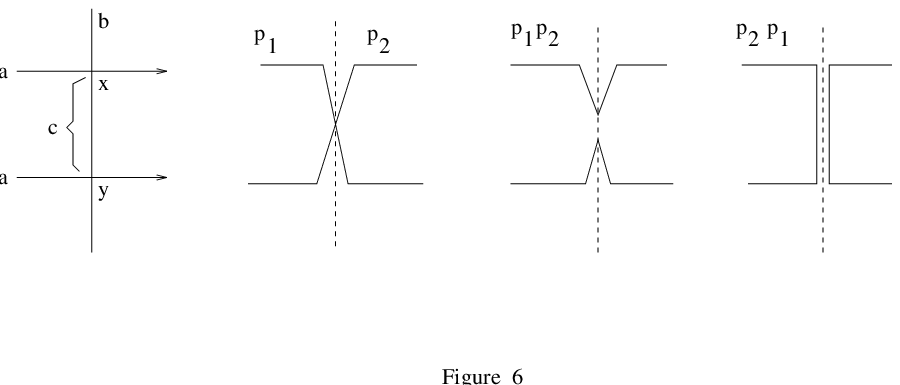}}
\midspace{0.1cm}
\centerline{\epsfbox{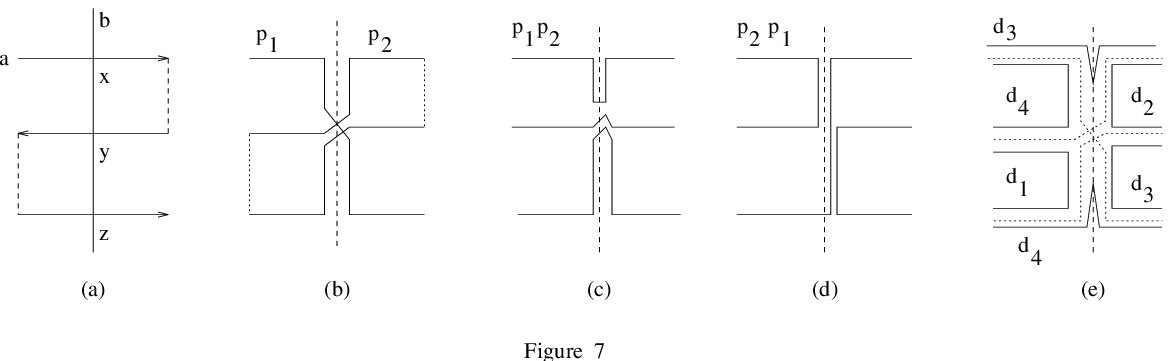}}
\midspace{0.1cm}

\it Case 2. \rm   Suppose any  pair of adjacent  intersection points 
in $b$ has different intersection signs. Then $|a \cap b| \geq 3$. 
Take three intersection
points $x_1, x_2, x_3$ in $b$ so that $x$, $y$ and $y$, $z$ are adjacent.
Their intersection signs alternate. Fix an orientation on $a$ so that
the arc from $x$ to $y$ in $a$ does not contain $z$ as shown in figure 7(a). 
If  the surface $\Sigma$  is right hand oriented as in figure 7(a), 
take  $p_1$ and $p_2$ as in figure 7(b). Then $p_1 \perp_0 p_2$
in $ N(p_1 \cup p_2)$. We claim that $p_1 \perp_0 p_2$ in
$\Sigma$. To see this, it suffices to show 
that $N(p_1 \cup p_2)$ is incompressible in $\Sigma$. Indeed, each
 boundary components of $N(p_1 \cup p_2)$ is  isotopic to 
a simple loop $b_i$ made
by the arcs with ends $x_i, x_{i+1}$ along $a$ and $d$. Since
$|a \cap d| = I(a, d)$, these loops $b_i$ are essential and non-boundary
parallel. Thus the claim  follows. By the construction
conditions (1), (2) and (3) follow from figure 7(c), (d) and (e).
If $\Sigma$ is left-hand oriented, we simply interchange $p_1$ and $p_2$.
$\square$

As an application of the proposition, we show that the mapping
class group is finitely generated by Dehn twists. Take $f$ in the
proposition to be the
map sending $\alpha \in \Cal S(\Sigma)$
to the isotopy class of positive Dehn twist along $\alpha$. 
First of all, there are two basic relations
on the Dehn twists: (1) (braid relation)
$D_{\alpha \beta } = D_{\alpha} D_{\beta} D_{\alpha}^{-1}$ for
$\alpha \perp  \beta$ and (2) (lantern relations) $ D_{\alpha} D_{\beta}
D_{\alpha \beta} = D_{\partial N(\alpha \cup \beta)}$ for
$\alpha \perp_0 \beta$. Thus by the proposition, the mapping 
class group is generated by elements in $G_0$.
For all surfaces, it is easy to construct a finite set $G_0$ so that
$G_{\infty} = \Cal S$. For instance, if
surface $\Sigma_{g,r}$ has $r >0$, let $
\{c_1,..., c_n \}$ ($n = 6g+3r-6$) be an ideal triangulation of it,
i.e, a maximal collection of disjoint pairwise non-isotopic, essential arcs in
 $\Sigma_{g,r}$. Then the corresponding collection $G_0$ in the
corollary is a finite set, indeed $|G_0|
\leq 3^n$ since each $\alpha \in CS(\Sigma)$ is determined by the
n-tuple $(I(\alpha, [c_1]),..., I(\alpha, [c_n])$). 

\it Remark. \rm The lantern relation was discovered  and used
by M. Dehn ([De], p333) and rediscovered independently by Johnson [Jo3].
Also the braid relation (1) implies the Artin's relation
$D_{\alpha}D_{\beta}D_{\alpha} = D_{\alpha}D_{\beta}D_{\alpha}$.

{\bf \S4. Thurston's Embedding of the Teichm\"uller Space}

We prove theorem 1  for compact surface $\Sigma_{g,r}$ in sections \S 4.1-4.3.
In \S 4.4, we indicate the modification needed for non-compact surfaces.
By the proof of theorem 2, it suffices to show
that conditions (1) and (2) are sufficient.

4.1. \it Reduction to the surfaces $\Sigma_{0,5}$ and $\Sigma_{1,2}$ \rm

We shall prove theorem 1 by  induction on the norm $|\Sigma_{g, r}|$ $
= 3g+r$ of a compact surface. The goal of this section is to
show that theorem 1 for all surfaces follows from theorem 1 for
$\Sigma_{0,5}$ and $\Sigma_{1,2}$.

Given $\Sigma =\Sigma_{g,r}$ with 
$|\Sigma| \geq 5$, and a function $f: \Cal S(\Sigma)
\to \bold R$ which is a trace function on each incompressible
subsurface $\Sigma'$ of norm 4, we decompose $\Sigma = X \cup Y$
so that $X$, $Y$ are incompressible of smaller norms with 
$ int(X \cap Y ) \cong int(\Sigma_{0,3})$ as figure 3(d). 
To be more precise, we take $X=\Sigma_{0,r-1}$, $Y= \Sigma_{0,4}$
if $g=0$ and take $X = \Sigma_{g-1, r+2}$, $Y = \Sigma_{1,1}$ if $g \geq 1$.
Consider the restrictions $f|_{\Cal S(X)}$ and $f|_{\Cal S(Y)}$. By the
induction hypothesis we find hyperbolic metrics $m_X$ and $m_Y$ on
$X$ and $Y$ respectively realizing the restrictions as the trace functions.
By the gluing lemma, we construct a hyperbolic metric $m$ on
$\Sigma$ whose restriction to $X$ and $Y$ are isotopic to $m_X$
and $m_Y$. Thus  the trace function $t_m$ and $f$ have the same
values on $\Cal S(X) \cup \Cal S(Y)$.

The goal is to show that  the above
condition  $f|_{\Cal S(X) \cup \Cal S(Y)} = t_m|_{\Cal S(X) \cup \Cal S(Y)}$
implies $f = t_m$.
To achieve this, let us rewrite the conditions (1), (2)
 satisfied by $f$ and $t_m$ as follows:

$ (1') \quad  f^2(\alpha)+ f^2(\beta) + f^2(\alpha \beta) - f(\alpha)f(\beta)
f(\alpha \beta) -2 + f(\partial N(\alpha \cup \beta)) = 0, \quad \quad $ if
 $\alpha \perp \beta$,

$ (2')  \quad f^2(\alpha) + f^2(\beta) +$$ f^2(\alpha \beta) -$
$ f(\alpha) f(\beta)$
$ f(\alpha \beta) + $$f(\alpha)(f(\gamma_1)f(\gamma_2) +$
$ f(\gamma_3) f(\gamma_4))+$$ f(\beta)$$(f(\gamma_2)$ \newline
$f(\gamma_3)$ + $ f(\gamma_1)f(\gamma_4)) +$$
f(\alpha \beta)(f(\gamma_2)f(\gamma_4)$$ + f(\gamma_1)f(\gamma_3))$$
+ f^2(\gamma_1)$$ + f^2(\gamma_2) +$$ f^2(\gamma_3)$$ + f^2(\gamma_4)$$
+ f(\gamma_1)$ \newline $f(\gamma_2)$ 
$f(\gamma_3)$$ f(\gamma_4) -4 = 0, $ if $\alpha \perp_0 \beta$,

$ (3') \quad f(\alpha \beta) + f(\beta \alpha) = f(\alpha) f(\beta), $ if $
\alpha \perp \beta,$ and

$ (4') \quad  f(\alpha \beta) + f(\beta \alpha) = f(\alpha)f(\beta) -
f(\gamma_1)f(\gamma_3) - f(\gamma_2)f(\gamma_4) $,  if $\alpha \perp_0 \beta,$

where  $\gamma_i's$ are the four components
of $\partial N(\alpha \cup \beta))$ so that $\alpha$ separates
\{$\gamma_1$, $\gamma_2$\} from \{$\gamma_3$, $\gamma_4$\}  and
$\beta$ separates $\{\gamma_2, \gamma_3\}$ and \{$\gamma_1, \gamma_4\}$.

Note that relations $(3')$ and $(4')$ give rise to an iteration
process. Namely, the value  $f(\alpha \beta)$ is
determined by the values of $f$ at $\alpha$, $\beta$, and
$\beta \alpha$ if $\alpha \perp \beta$, and is determined by
the values of $f$ at $\alpha$, $\beta$, $\beta \alpha$ and
the four components of $\partial N(\alpha \cup \beta)$ if
$\alpha \perp_0 \beta$.

Let $a_1$, $a_2$ be the simple loops in $\partial (X \cap Y)$ which is
non boundary parallel in $\Sigma$ as in figure 3(d).
Applying proposition 1 to  $f$ and  to $t_m$ with respect to the set
$\{a_2\}$, we conclude that 
$f$  = $t_m$  follows from $f(\alpha) =t_m(\alpha)$ where
$\alpha \perp_0 [a_2]$. 
Assume that theorem 1 holds for $\Sigma_{0,5}, \Sigma_{1,2}$. We show 
$f(\alpha) = t_m(\alpha)$ with $\alpha \perp_0 [a_2]$ as follows.
Take $s \in \alpha$ so that $|s \cap a_2|=2$. Then $Z = Y \cup N(s)$
is an incompressible 
subsurface homeomorphic either to $\Sigma_{1,2}$ or $\Sigma_{0,5}$.
Let $X' = X \cap Z$, $Y' = Y \cap Z$. Then $Z = X' \cup Y'$ so that
$X' \cap Y' = X \cap Y$. Consider $f|_{\Cal S(Z)}$ and $t_m|_{\Cal S(Z)}$.
By theorem 1 for $Z$ and the fact that $f$ and $t_m$ coincide on
the subset $\Cal S(X') \cup \Cal S(Y')$, we conclude that
$f = t_m$ on $\Cal S(Z)$ by the gluing lemma. In
particular, $f(\alpha) = t_m(\alpha)$.

It remains to show theorem 1 for $\Sigma_{0,5}$ and $\Sigma_{0,5}$.
By the same decomposition $\Sigma = X \cup Y$ as above, it suffices
to show the following two lemmas.

For simplicity, we let $Im(\Sigma)$ be the set of all functions
from $\Cal S(\Sigma)$ to $\bold R_{>2}$ satisfying conditions $(1')$, $(2')$, $(3')$,
and $(4')$. Two classes $\alpha$ and $\beta$ are \it disjoint \rm if they
are distinct and have disjoint representatives.

{\bf Lemma 8.} \it Suppose $\alpha_1$ and $\alpha_2$ are two 
 disjoint elements in $\Cal S'(\Sigma_{0,5})$. If two
elements $f$ and $g$ in $Im(\Sigma_{0,5})$ satisfy $f(\alpha) = g(\alpha)$
for all $\alpha \in \Cal S(\Sigma_{0,5})$ with $I(\alpha, \alpha_1)
I(\alpha, \alpha_2) = 0$, then $f =g$. \rm

{\bf Lemma 9.} \it Suppose $\alpha_1$ and $\alpha_2$ are two
disjoint  elements in $\Cal S'(\Sigma_{1,2})$ so that
$\alpha_1$ is non-separating and $\alpha_2$ is separating. If $f$ and $g$
are two elements in $Im(\Sigma_{1,2})$ so that $f(\alpha) = g(\alpha)$
for all $\alpha \in \Cal S(\Sigma_{1,2})$ with $I(\alpha, \alpha_1)
I(\alpha, \alpha_2) = 0$, then $f = g$.
\rm

4.2. \it Proof of lemma 8 \rm

To prove lemma 8, by  proposition 1,
it suffices to show
that $f(\alpha) = g(\alpha)$ for $\alpha \perp_0  \alpha_i$ for $i = 1,2$.
Let $a_i \in \alpha_i$ be a representative so that $|a_1 \cap a_2| = 0$ 
and let
$x \in \alpha$ so that $x \perp_0 a_i$ for $i =1,2$. Note that if
$x' \perp_0 a_i$ for $i =1,2$, there is an orientation preserving
homeomorphism $h$ of $\Sigma_{0,5}$ sending $x$ to $x'$ and
preserving each $a_i$ (since both
$N(a_1 \cup a_2 \cup x)$ and $N(a_1 \cup a_2 \cup x')$ are strong
deformation retractors for $\Sigma_{0,5}$). Thus we may draw $x$
as in figure 8(a). Let $a,b,c,d,e$ and $b_1$, $b_2$, $c_1$, $c_2$, $d_1$,
and $d_2$ be curves as in figures 8(a), (b) and (c) 
so that each of them is either  disjoint from $a_1$ or from $a_2$.

\midspace{0.1cm}
\centerline{\epsfbox{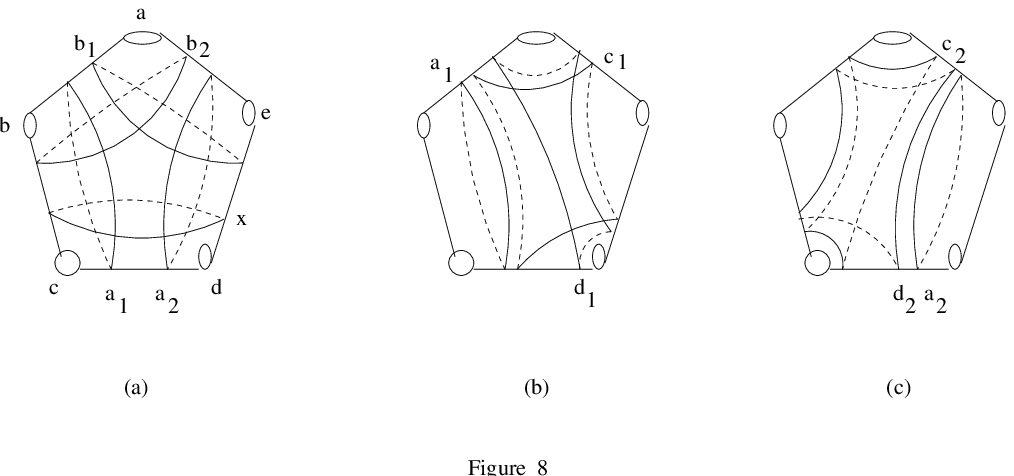}}
\centerline{\epsfbox{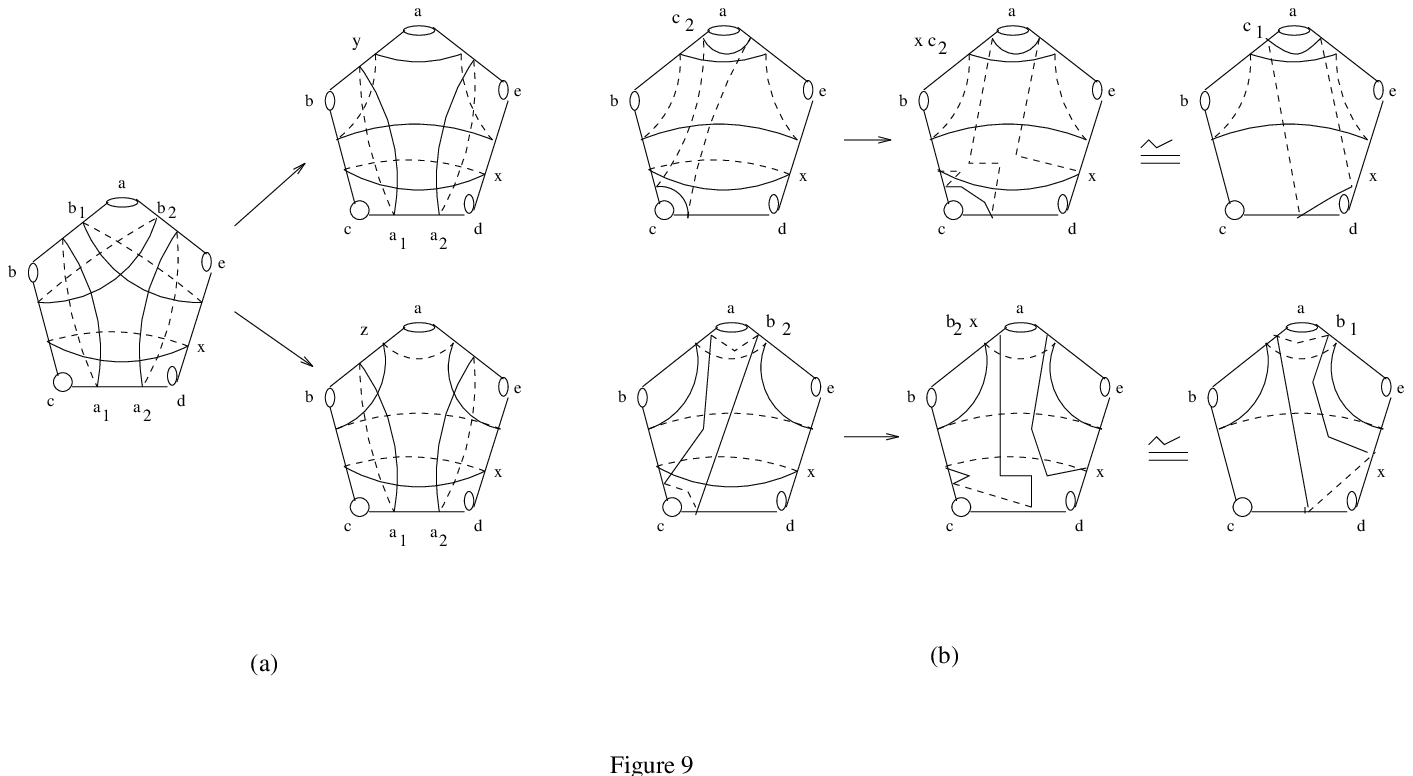}}
\midspace{0.1cm}

{\bf Claim.} \it There is a rational function  $R$ so that
for each $h \in Im(\Sigma_{0,5}$), $h(x )  = R(h(a), h(b),$ ..., $h(e)$,
$h(a_1), h(a_2),$ $ h(b_1), h(b_2), h(c_1), h(c_2), h(d_1), h(d_2))$.
\rm

It follows from the claim that $f(\alpha) = g(\alpha)$.
This finishes the proof of lemma 8.

Before begin the proof of the claim, let us simplify the notations by making the
following conventions.

(C1) The value of $h$ at a curve $s$ will be denoted by  $s$.

(C2) The multiplication of two curves $s_1$ and $s_2$ will be denoted by $s_1
\circ s_2$.

(C3) Surfaces drawn in the figures have the right-hand orientation in the
front face.

Let $y =b_1 \circ b_2$ and $z = b_2 \circ b_1$ as in figure 9(a).

Since $b_1 \perp_0 b_2$ and  $\partial N(b_1 \cup b_2) \cong a \cup b \cup e
 \cup x$, applying relation ($2'$) in $N(b_1 \cup b_2)$ with 
respect to $b_1$, $b_2$ and $y$,
we obtain:
$ x^2 + a^2 + b^2 + e^2 + y^2 + b_1^2 + b_2^2 - b_1b_2y + abex
+b_1(ae + bx) + b_2(ex + ab) + y(ax + be) -4 = 0.$
This can be written as:
$$ x^2 + y^2 + axy + p_1x - p_2y + p_3 = 0,    \tag 15$$
where $p_j$ are some polynomials in $a,b,c,d,e, a_i, b_i, c_i,$ and
$d_i$ (the same notations apply below) and $p_j >0$ for $j = 1,2,3$.

Similarly, we have,
$$ x^2 + z^2 + axz + p_1x - p_2z + p_3 = 0.  \tag 16$$
Furthermore, by ($4'$), $y+z = b_1b_2 - ax - be$, i.e.,
$$ax + y + z = p_4.   \tag 17$$
Now $c_2 \perp_0 x$ and $x \circ c_2 = c_1$ (see figure 9(b)). Applying the
relation ($2'$) to $N(c_2 \cup x)$ with respect to
$c_2, x, c_1$ and using $\partial N(c_2 \cup x) \cong a \cup c \cup d
\cup y$,  we obtain $y^2 + a^2 + c^2 + d^2 + x^2 + c_1^2 + c_2^2
-c_1 c_2 x + acdy + x(ay + cd) + c_2(ac + dy) + c_1(ad + cy) -4 = 0$, i.e.,
$$x^2 + y^2 + axy -p_5x + p_6y + p_7 = 0,  \tag 18$$
where $p_5$, $p_6$ and $p_7$ are positive.

Similarly, using $d_2 \perp_0 x$ and $d_1 \circ x = d_2$, we obtain 
a relation:
$$x^2 + z^2 + axz -p_8x + p_9z + p_{10} = 0,  \tag 19$$
where $p_8$, $p_9$ and $p_{10}$ are positive.

Consider the difference of (15) and (18). We obtain,
$$ p_{11}  x - p_{12} y = p_{13}, \tag 20$$
where $p_{11}$ and $p_{12}$ are positive.

Consider (16)-(19), we obtain,
$$p_{14} x - p_{15} z = p_{16},  \tag 21$$
where $p_{14}$ and $p_{15}$ are positive. 

Now the system of linear equations (17), (20) and (21) in variables
$x, y$ and $z$ has a unique solution since its determinant is positive.
This ends to the proof of the claim and thus finishes the proof of lemma 8.

4.3.  \it Proof of lemma 9 \rm

To prove lemma 9, by proposition 1,
it suffices to
show that $f(\alpha) = g(\alpha)$ for $\alpha \in \Cal S(\Sigma_{1,2})$ with
$\alpha \perp_0 \alpha_2 $ and $
\alpha \perp \alpha_1 $ since there is no element $\beta \in \Cal S(\Sigma_{1,2})$
so that $\beta \perp_0 \alpha_i $ for $i =1,2$. Fix such an $\alpha$
for the rest of the proof. Take $x \in \alpha$,
$a_i \in \alpha_i$, $i=1,2$ so that $a_1 \cap a_2 = \emptyset$,
$x \perp a_1$ and  $x \perp_0 a_2$. 

Let $Y = \Sigma_{1,2} - a_1$ and $X$ is the subsurface bounded
by $a_2$ containing $a_1$. We  have $f = g$ on the subset $\Cal S(X) \cup
\Cal S(Y)$.

{\bf Claim.}  \it There exists  a finit  set of elements 
$\{\beta_1,..., \beta_n\}$ in $\Cal S(X) \cup
\Cal S(Y)$ and a function $F$ so that
for any element  $h$ in $Im(\Sigma_{1,2})$, $h(\alpha)$ = $F( h(\beta_1),...,
h(\beta_n))$.  \rm 

It follows  from the claim that $f(\alpha) = g(\alpha)$. This finishes the
proof of lemma 9.

We shall adopt the same
convention as in \S4.2 by identifying $h(s)$ with  the simple closed curve
$s$ for the rest of the proof.

{\bf Proof of the claim}.
Since any other simple closed curve $x'$  with $x' \perp a_1$ and $x' \perp_0
 a_2$
is an image of $x$ under an orientation preserving self-homeomorphism
preserving $a_1$ and $a_2$, we may draw $x$ as in figure 10. Introduce 
a few more curves $y$, $z$,  $x_1$, $y_1$, $x_2$, $y_2$, 
$b_1$, $b_2$, $b_3$, $k$ as in figure 10. Note that
the curves $b_1$, $b_2$, $b_3$,  and $k$
are either in $X$ or in $Y$.

\midspace{0.1cm}
\centerline{\epsfbox{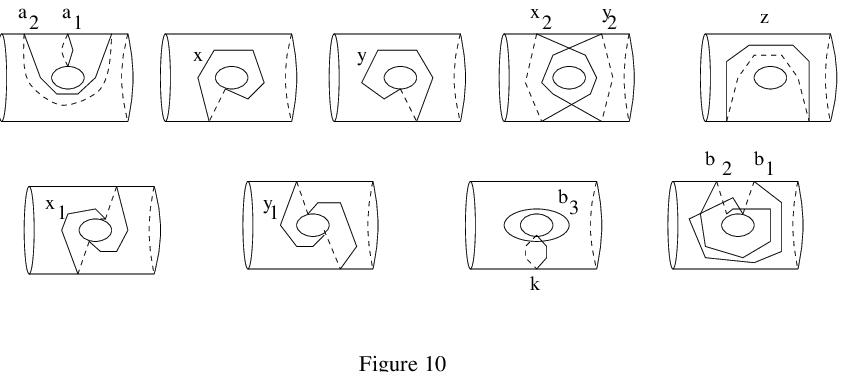}}
\centerline{\epsfbox{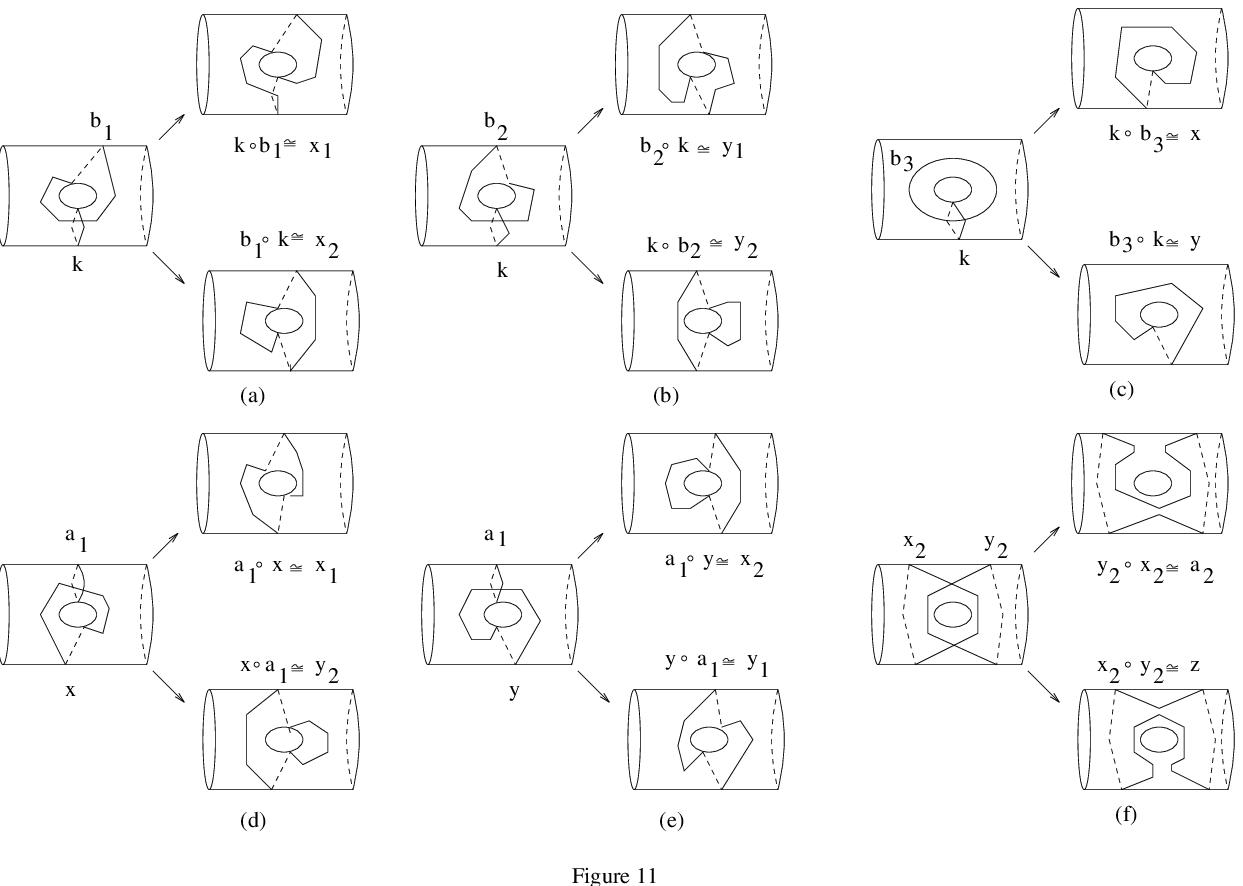}}
\midspace{0.1cm}

There are many relations among these curves as shown in figure 11.

We obtain a system of equations in $x, y, x_1, y_1, x_2, y_2$ and $z$ by
applying  formulas $(1')$, $(2')$, $(3')$ and $(4')$.

By figure 11(a), we have,
$$ x_1 + x_2 = kb_1.  \tag e1$$
By figure 11(b), we have,
$$y_1 + y_2 = kb_2.  \tag e2$$
By figure 11(c), we have,
$$ x + y = kb_3.  \tag e3$$
By figure 11(d) and that $\partial N(x \cup a_1) \subset Y$, we have,
$$  x_1 + y_2 = a_1 x, \tag e4$$
and,
$$ y_2^2 + x^2 - a_1 xy_2 = p_1.    \tag e5$$

By figure 11(e) and that $\partial N( y \cup a_1) \subset Y$, we have,
$$ y_1^2 + y^2 - a_1 y y_1 = p_2.  \tag e6 $$
By figure 11(c) that $x = k \circ b_3$ and 
$\partial N(b_3 \cup k) \cong z$,  we have,
$$ x^2 - kb_3x = -z + p_3. \tag e7$$
By figure 11(f) and $\partial N(x_2 \cup y_2) \cong a \cup b \cup b_3 \cup b_3'$
where $b_3'$ is a parallel copy of $b_3$ and $\partial \Sigma_{1,2} 
= a \cup b$, we have,
$$ z = x_2y_2 + p_4. \tag e8$$
Here and below, $p_i's$  denote 
some polynomials in  some elements in $\Cal S(X) \cup \Cal S(Y)$.

Also from $ a_1 \perp b_3$ with $b_3 \circ a_1 = b_2$ and 
$a_1 \circ b_3 = b_1$, we have,
$$ b_1 + b_2 = a_1 b_3. \tag e9$$

The goal is to show that the system of equations (e1)-(e8) has 
a unique solution in $x$. Assuming
this, we conclude that  the claim holds.

To this end, we shall first eliminate $x_1$, $y_1$,
$y$ and $z$ from the above system and show that $x_2$ are $y_2$ are
linear functions in $x$.

Subtracting (e1) by (e4) gives:
$$x_2 = y_2 - a_1 x + kb_1,  \tag e10$$
and subtracting (e7) by (e8) gives:
$$x_2y_2 + x^2 -kb_3x = p_5.   \tag e11$$
By  (e3), $y = kb_3 - x$ and by (e2), $y_1 = kb_2 - y_2$. Substitute them into
(e6) and subtract the result by (e5), we obtain:
$$ (a_1 b_3 - 2 b_2) ky_2 + (a_1b_2 - 2b_3) kx = p_6. \tag e12$$
Note that the coefficients of $y_2$ and $x$ in (e12) cannot be both zero since
$a_1 > 2$. If $a_1b_3 - 2b_2 =0$, then $x$ is determined uniquely. Suppose
otherwise, then we solve $y_2$ in terms of $x$ and obtain,
$$ y_2 =  p_7 x +  p_8.  \tag e13$$
where $p_7 = (2b_3 -a_1b_2)/(a_1 b_3 - 2b_2)$.
From (e10), we obtain
$$ x_2 = ( p_7 - a_1 ) x +  p_8 + kb_1.  \tag e14$$
Now substitute (e13) into (e5), we obtain a quadratic equation in $x$ as
follows:
$$ (p_7^2 - a_1 p_7 + 1) x^2 + ( 2 p_7 p_8 -a_1 p_8) x + p_9 = 0.  \tag e15$$

Substitute $(e13)$ and $(e14)$ into $(e11)$ to obtain a quadratic
equation  in $x$ as follows.

$$ (p_7^2 -a_1 p_7 + 1) x^2 + (-b_3k + p_8(p_7 -a_1) + p_7p_8 + kb_1 p_7)x
+ p_{10} = 0. \tag e16$$

Subtract (e16) by (e15) to obtain a linear equation in $x$ whose
leading term is $-k b_3 + k b_1 p_7$.
Replace $p_7$ by $(2b_3 -a_1b_2)/(a_1 b_3 - 2b_2)$ and use (e9)
that $b_1 = a_1b_3 - b_2$, we simplify the leading coefficient to
$a_1k(b_2^2 + b_3^2 - a_1b_2b_3)/(a_1b_3 - 2b_2)$. The number
$b_2^2 + b_3^2 - a_1b_2b_3$ is negative  by relation $(1')$  that
$a_1^2 + b_2^2 + b_3^2 < a_1 b_2 b_3$. Thus we obtain a unique
solution of $x$.
This finishs the proof of lemma 9. 
$\square$

4.4.  \it Proof of  Theorem 1 for metrics with cups ends \rm

We first recall theorem 2 for metrics with cups ends. Let $\Sigma =
\Sigma_{0,r}^s$
with $r+s =4, s<4$, be given with three simple closed curves
$a_{12}, a_{23},$ and $a_{31}$ on it satisfying $a_{31} = a_{12}a_{23}$ and
$a_{12} \perp_0 a_{23}$. Let $b_i$ be four essential simple closed curves in
$int(\Sigma_{0,r}^s)$ which are homotopic into the four ends so that
$a_{ij}$, $b_i$ and  $b_j$ bound a 3-holed sphere in the surface ($i \neq j,
i,j \leq 3$).
Assume the cusp ends correspond to $b_i$ ($i=1,2,.., s$).
Take the collection $F \subset \Cal S(\Sigma)$  
to be the isotopy classes of $a_{ij}$ and $b_i's$ where $i \neq j$ and
$i,j \leq 3$. Then   the same
argument used in the proof of theorem 2 shows,

{\bf Lemma 10.} \it  The map $\pi_F : T(\Sigma_{0,r}^s) \to  \bold R_{\geq 2}^6$is an embedding whose image is given by $\{(t_1, t_2, t_3, $$t_{12},$$
t_{23}, t_{31})$$\in \bold R_{\geq 2}^6 |$ $t_1 = ... =t_s = 2$,
 $t_{s+1} >2,..., t_3 > 2$, so that the formula (4) holds\}. Furthermore,
 there exists a real analytic RC-section for $T(\Sigma_{0,r}^s)$. \rm 

Now to construct metrics on $\Sigma_{g,r}^s$ with $s>0$, 
we use the decomposition $\Sigma_{g,r}^s = X \cup Y$ as in figure 12.
 The first case  (1) is given by $r>0$.  We need to consider subcases
(1.1), (1.2) and (1.3) where (1.1) corresponds to $g>0$, (1.2)
corresponds to $g=0$ and $r+s>5$, and (1.3) corresponds to
$g=0$ and $r+s \leq 5$.   In cases (1.1), or  (1.2), we
choose $X \cong \Sigma_{g, r}^{s-1}$, $Y \cong \Sigma_{0,3}^1$, and
$X \cap Y \cong \Sigma_{0,3}$.  In case (1.3) with $r+s \leq 4$, then
it follows from theorem 2. In case (1.3) and $r+s =5$,
we choose $X \cong \Sigma_{0, r-1}^{s-1}$, $Y \cong \Sigma_{0,u}^v$,
where $u+v = 4$, $ 2 \geq v \geq 1$, and $X \cap Y \cong \Sigma_{0, 4-v}^{v-1}$.
In the second case (2) $r=0$, we need to consider subcases 
(2.1) $ s \geq 2$ and (2.2) $s =1$.  
In case (2.1) that $s \geq 2$, if (2.1.1) $g >0$, or (2.1.2)
$ g =0$ and $s >5$, then $X \cong \Sigma_{g,1}^{s-2}$, 
$Y \cong \Sigma_{0,2}^2$, and $X \cap Y \cong \Sigma_{0,3}$. 
If  (2.1.3) $ 5 \geq s \geq 2$ and $g=0$, 
the theorem holds except for $s=5$ where
we decompose $\Sigma_{0,0}^5$ as a union of two $\Sigma_{0,1}^3$ with
intersection $\Sigma_{0,2}^1$.  Finally, in case (2.2) that $s=1$, it suffices
to consider $g \geq 2$. 
We take $X \cong \Sigma_{1, 1}$, $Y \cong \Sigma_{g,1}^1$
and $X \cap Y \cong X -s$ where $s$ is a non-separating simple closed curve
in $X$.

\midspace{0.1cm}
\centerline{\epsfbox{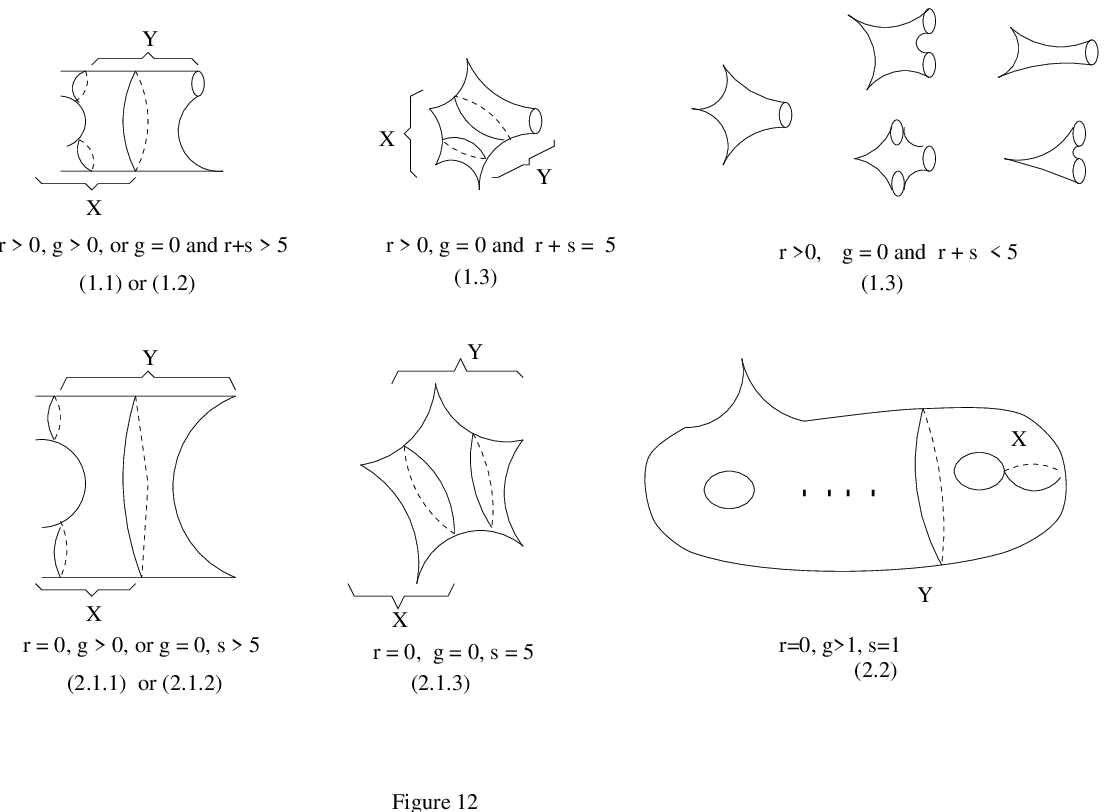}}
\midspace{0.1cm}

These give the 3-holed  sphere decomposition of the surface into 
two subsurfaces of smaller $|X|$ and $|Y|$ where $|\Sigma_{g,r}^s| = 3g+r+s$.
Note that lemmas 8 and 9 still hold for metrics with cups ends.
Now by the gluing lemma, lemmas 8, 9, 10, theorem 2, the same 
argument used in the previous sections applies. This gives a proof of theorem 1 
for metrics with cusp ends.

\it Remark.  \rm Teichm\"uller space is well known to be homeomorphic to a Euclidean space.
This fact can also be  derived from theorem 2  and lemma 1.
Indeed, the gluing lemma  shows that the restriction map from
 $T(X \cup Y)$ to $T(X)$ is a fiber-bundle map.
The fiber can be  shown to be homeomorphic to a Euclidean space by solving
a simple inequality (e.g. relations (3) or (4)).

{\bf  \S5. Applicaion to Finite Dimensional Embeddings of Teichm\"uller
Spaces}

We shall prove  the following stronger version of the corollary
for compact surfaces by 
induction on $|\Sigma_{g,r}| = 3g+r$ in this section. The proof for
surfaces with cusp ends will be omitted. 

{\bf Corollary.} \it (a)
For  surface $\Sigma_{g,r}$ of negative Euler number and $r >0$, there
exists a finite subset $F$ in $\Cal S(\Sigma_{g,r})$ consisting of
$6g + 3r -6$ elements so that the map $\pi_F : T(\Sigma_{g,r})
\to \bold R_{>2}^F$ is an embedding onto an open subset which is
defined by a finite set of real analytic RC-inequalities in the
coordinates of $\pi_F$. Furthermore, there exists a map $f: T(\Sigma_{g,r})
\to \tilde R(\Sigma_{g,r})$ so that for each $m$ in $T(\Sigma_{g,r})$,
$f(m)$ is a lifting of a monodromy of $m$ and  the entries of
the matrix $f(m)(\alpha)$  are real analytic RC-functions of $\pi_F(m)$ for any 
$\alpha \in \Cal S(\Sigma_{g,r})$.

(b) For surface $\Sigma_{g, 0}$ of negative Euler number, there exists
a finite subset $F$ of $\Cal S(\Sigma_{g,0})$ consisting of $6g-5 $
elements so that $\pi_F : T(\Sigma_{g,0}) \to \bold R_{>2}^F$ is an 
embedding whose image is defined
by one real analytic  RC-equation and finitely many real analytic RC-inequalities 
in the coordinates of $\tau_F$.
Furthermore, there exists a map $f: T(\Sigma_{g,0})
\to \tilde R(\Sigma_{g,0})$ so that for each $m$ in $T(\Sigma_{g,0})$,
$f(m)$ is a lifting of a monodromy of $m$ and  the entries of
the matrix $f(m)(\alpha)$  are real analytic RC-functions of $\pi_F(m)$ for any 
$\alpha \in \Cal S(\Sigma_{g,0})$
\rm

Note that the corollary without the statement about the lifting of monodromies
follows immediately from the gluing lemma, theorems 1 and 2, and lemmas 8 and 9.
To prove the full statement, we need to strengthen the gluing lemma.

In \S 5.1, we prove an extended version of the gluing lemma. In
\S 5.2, we prove the corollary  for $\Sigma_{1,2}$. The corollary for surfaces
with non-empty boundary is proved in \S 5.3. In \S 5.4, we prove the corollary
for closed surfaces.

5.1. \it Algebraic dependence in the gluing lemma \rm

We begin with a parametrized version of the Jordan
canonical form theorem for SL(2,$\bold R$) matrices.

{\bf Lemma 11.} \it (a) If A =$[a_{ij}]$ in SL(2,$\bold R$)  satisfies
$|tr A| > 2$ and $a_{12}a_{21} \neq 0$, then
$$C^{-1}AC =
1/2 \left(\matrix a_{11} + a_{22} + \sqrt{(a_{11} + a_{22})^2-4} & 0\\0& 
a_{11} + a_{22} - \sqrt{(a_{11} + a_{22})^2-4} \endmatrix \right),$$

where $$C =\left(\matrix  2a_{12} &  a_{11} -a_{22} - \sqrt{(a_{11} + a_{22})^2 -4}   \\
a_{22} -a_{11} + \sqrt{(a_{11} + a_{22})^2 -4} &  2a_{21} \endmatrix \right).$$

(b) For $A$ = $[a_{ij}]$ and $B = [b_{ij}]$ in SL(2,$\bold R$) with $trA > 2$
(resp. $trA < -2$), $a_{12} a_{21} \neq 0$ and $trABA^{-1} B^{-1} \neq 2$, 
there exist four real analytic RC-functions $c_{ij}$ in eight variables so that $C^{-1}AC =$
$\left(\matrix  \lambda & 0\\0&  \lambda^{-1}\endmatrix \right)$, $\lambda > 1$ (resp. $\lambda < -1$)
and $C^{-1}BC =$ $\left(\matrix  \alpha &  \beta\\\delta& 
\gamma\endmatrix \right)$, $|\delta| =1$ where $C = [c_{ij}(A,B)] \in$ GL(2,$\bold R)$.
\rm

{\bf Proof.} Part (a) follows by a direct calculation. Note that the
matrix $C$ is invertible since  $a_{12} a_{21} \neq 0$. 
Part (b) follows
from part (a). Indeed, by part (a), we may conjugate $A$ to the required
diagonal form $A'$. We also conjugate $B$ by the same matrix to obtain $B'$.
The trace of
the commutator remains unchanged. Thus the new matrix $B'$ = $[b'_{ij}]$
has non-zero (2,1)-entry. Now a further conjugation by  the
matrix $\left(\matrix  \sqrt{|b'_{21}|} & 0 \\ 0&  \sqrt{|b' _{21}|}^{-1}
\endmatrix \right)$ will
not change matrix $A'$ but change $B'$ into the required form.
$\square$

We say a pair of matrices $(A,B)$ is  \it normalized \rm
if $A = $   $\left(\matrix  \lambda & 0\\0&  \lambda^{-1}\endmatrix \right)$
with  $ |\lambda| > 1 $
and the (2,1)-entry of $B$  is 1. It follows from the
normalized condition that  if $C$ is in GL(2,$\bold R)$ so that
both $(A,B)$ and $(C^{-1}AC, C^{-1}BC)$ are normalized, then 
$(A,B) =$$(C^{-1}AC, C^{-1}BC)$, i.e.,  normalization is unique up to
GL(2, $\bold R$) conjugation.
Fix a pair  of elements $(\gamma_1, \gamma_2)$ in $\pi_1(\Sigma)$. A
representation $\rho$ in $\tilde R(\Sigma)$ is called \it normalized \rm with
respect to the pair if $(\rho(\gamma_1), \rho(\gamma_2))$ is normalized.
 
A  \it section \rm  of the natural projection from $\tilde R(\Sigma)$ to
$T(\Sigma)$ is a continuous map $f: T(\Sigma) \to \tilde R(\Sigma)$
so that $f(m)$ is a lifting of a monodromy of $m$.
Given a section $f$, we may produce a new section whose image
lies in any given component of $\tilde R(\Sigma)$ as follows. Conjugating
representations in $f(T(\Sigma))$ by the matrix
$\left(\matrix  1   &  0  \\ 0  &  -1 \endmatrix \right)$ 
gives rise to a new section in a component of the opposite orientation type;
and choosing a different lifting $\rho_{I}$ associated
to $\rho \in f(T(\Sigma))$
for a fixed index set $I$ (see \S 1.2 for the definition) gives
a section in a different component of the same orientation type.
We call these new sections  to be the ones obtained from 
$f$ by different liftings and
conjugations.
An  \it RC-section \rm
is a section so that (1) there exists an associated finite set $F$ $\subset 
\Cal S(\Sigma)$ so that the entries of the matrix $f(m)(\alpha)$ are
real analytic RC-functions
of the coordinates of $\pi_F(m)$ for all $\alpha \in \pi_1(\Sigma)$ and (2) each representation  in
the image of the section is normalized
with respect to a fixed pair of elements in $\pi_1(\Sigma)$. By lemmas
4, 5 and theorem 2, the Teichm\"uller spaces $T_{0,3}$, $T_{1,1}$, and
$T_{0,4}$ have RC-sections.

For simplicity, we shall identify curves, isotopy classes of curves,
and homotopy classes of curves in incompressible subsurfaces with their
images in the ambient spaces without mentioning the including maps.

{\bf Lemma 12.} (Algebraic dependence) \it Let $X$ and $Y$ be good
incompressible
subsurfaces of $\Sigma$ so that $\Sigma = X \cup Y$ and either (1) 
$X \cap Y \cong \Sigma_{0,3}$, or (2) $Y \cong \Sigma_{1,1}$ and 
$X \cap Y$ = $Y -s$
where $s$ is a non-separating simple closed curve in int($Y$),
or (3) $X \cap Y \cong \Sigma_{0,2}^1$ so that the punctured end
in $\Sigma_{0,2}^1$ is a punctured end of $\Sigma$. If
$T(X)$ and $ T(Y) $ both have RC-sections $f_X$ and $f_Y$ with associated
sets $F_X$ and $F_Y$ respectively, then $T(\Sigma)$ has an RC-section
with associated set $F_X \cup F_Y$. \rm

{\bf Proof.} Let $(\alpha_1, \alpha_2)$ (resp. $(\beta_1, \beta_2)$) be the pair in $\pi_1(X)$
(resp. $\pi_1(Y)$) so that each representation in the image of $f_X$
(resp. $f_Y$) is normalized with respect to it.
Choose two
geometric generators $\gamma_1$ and $\gamma_2$ for $\pi_1(X \cap Y)$ 
so that $\gamma_1 \gamma_2$
is represented by the third boundary component. Then one of the
three elements  $\gamma_1$, $\gamma_2$, $\gamma_1 \gamma_2$, say
$\gamma_1$, satisfies the condition that
both subgroups $<\alpha_1, \gamma_1>$ and $< \beta_1, \gamma_1>$ are not solvable. 
Let $\gamma_2$ be one of the remaining element. Then
$\pi_1(X \cap Y)$  is generated by  $\gamma_1$ and  $\gamma_2$.
We extend \{$\gamma_1$, $\gamma_2$\} to a minimal set of generators
\{$\gamma_1$, ..., $\gamma_n$\} for $\pi_1(\Sigma)$ so that
each $\gamma_i$ is either in $\pi_1(X)$ or in $\pi_1(Y)$.

By choosing a different lifting if necessary, we may assume that 
 $f_X (m)(\alpha_1)$ and $f_Y(m)(\beta_1)$
are diagonal  matrices with positive traces for $m \in T(X \cap Y)$
($f_X$ and $f_Y$ are still sections but may not be normalized any more).
Now by the choice of element $\gamma_1$, 
both matrices $f_X(m)(\gamma_1)$ and
$f_Y(m)(\gamma_1)$ have non-zero  off diagonal entries for all $m$, and
the trace of the commutator of $f_X(m)(\gamma_1)$ and $f_X(m)(\gamma_2)$
(resp.  $f_Y(m)(\gamma_1)$, and $f_Y(m)(\gamma_2)$) is not $2$. Thus
by lemma  11, we may conjugate
the pair $( f_X(m)(\gamma_1), f_X(m)(\gamma_2) )$ (resp. $(f_Y(m)(\gamma_1),
f_Y(m)(\gamma_2)))$  to the form in lemma 11(b) by a GL(2, $\bold R)$ 
matrix whose entries are real analytic
RC-functions in  the coordinates of $\pi_{F_X}(m)$ 
(resp. in the coordinates of $\pi_{F_Y}(m)$).
This produces two sections $g_X$ and $g_Y$ for $T(X)$ and $T(Y)$ respectively
so that (1) for each $m \in T(X)$ (resp. $m \in T(Y)$), the entries of the
matrices $g_X(m)(\gamma)$  (resp. $g_Y(m)(\gamma)$) are real analytic
RC-functions 
in the coordinates of $\pi_{F_X}(m)$ (resp. $\pi_{F_Y}(m)$),
and (2) the matrix $g_X(m)(\gamma_1)$  (resp. $g_Y(m)(\gamma_1))$
is diagonal with 
(1,1)-entry bigger than one and the (2,1)-entry of $g_X(m)(\gamma_2)$ 
(resp. $g_Y(m)(\gamma_2)$) has absolute value one.

We may normalize the sections $g_X$ and $g_Y$  by choosing different lifting
 which changes the generator $\rho(\gamma_2)$ to $-\rho(\gamma_2)$ if necessary.
Thus we may assume that both $g_X$ and $g_Y$ are normalized with respect to
the pair $(\gamma_1, \gamma_2)$. 

We now define an RC-section for $T(\Sigma)$ as follows. By the gluing lemma,
each $m \in T(\Sigma)$ corresponds to a pair $(m_X, m_Y) \in T(X) \times
T(Y)$ so that $R_{X }(m) = m_X$,  $R_{Y} (m) = m_Y$  and the
restrictions of $m_X$ and $m_Y$ to $X \cap Y$ are the same.
The  restrictions of the two representations $g_X(m_X)$ and $g_Y(m_Y)$ 
to the subgroup $\pi_1(X \cap Y)$ uniformize the same element $R_{X \cap Y}(m)$.
Since the pair $(\gamma_1, \gamma_2)$ generates $\pi_1(X \cap Y)$,
by the normalization condition for $g_X$ and $g_Y$, we have
$g_X(m_X)|_{\pi_1(X \cap Y)}$  $= g_Y(m_Y)|_{\pi_1(X \cap Y)}$. 
By Maskit combination theorem (there is no need to
verify the side condition since the gluing is along a 3-holed sphere),
there exists a unique representation $\rho \in \tilde R(\Sigma)$
so that $\rho |_{\pi_1(X)}$
$= g_X(m)$ and $\rho |_{\pi_1(Y})$ $ = g_Y(m)$.  The map from
$T(\Sigma)$ to $\tilde R(\Sigma)$ sending $m$ to $\rho$
is a section normalized with respect
to $(\gamma_1, \gamma_2)$. To see the RC-dependence (which also shows
the continuity of the map $m$ to $\rho$), it suffices to
check the condition for each generator $\gamma_i$. By
the  construction, $\rho(\gamma_i)$
is either $g_X(m)(\gamma_i)$ or $g_Y(m)(\gamma_i)$. Thus, each
entry of the matrix $\rho(\gamma_i)$ is a real analytic
RC-function in the coordinates
of $\pi_{F_X \cup F_Y}(m)$. 
$\square$

5.2. \it Proof of the Corollary  for $ \Sigma_{1,2}$ \rm

Let $s_7$ be an essential separating simple closed curve
and $s_1$ be a non-separating simple closed curve disjoint from
$s_7$ in $\Sigma_{1,2}$ as in figure 13(a). We decompose $\Sigma_{1,2}$ as
a union $X \cup Y$ where $X$ is the compact subsurface bounded by $s_7$
containing $s_1$ and $Y$ is the complement of $s_1$. Then $X \cap Y$
is $X - s_1$. Let $s_2$, $s_3$ be simple closed curves in
$X$ so that $s_1 \perp s_2$ and $s_3 = s_1 s_2$; let $s_4$, $s_5$,
$s_6$, $s_1^+$ and $s_1^-$ be simple closed curves in $Y$ so that
$s_1^+$ and  $s_1^-$ are boundary components which are identified to be
$s_1$ in $\Sigma_{1,2}$, $s_6  \subset \partial Y$ and $s_4 \perp_0 s_7$, $s_5 = s_4 s_7$.
See figures 13(b) and (c).  By the gluing lemma and lemma 4, the
Teichm\"uller space $T(\Sigma_{1,2})$ can be identified with the
subset $\{(m_X, m_Y) \in T(X) \times T(Y) | t_{m_X}(s_1) = t_{m_Y}(s_1^+)
=t_{m_Y}(s_1^-)$ and $t_{m_X}(s_7) = t_{m_Y}(s_7)$\}. 
By lemma 5,
$m_X$ is determined by $\pi_{F_X}(m_X)$$= (t_{m_X}(s_1),
t_{m_X}(s_2), t_{m_X}(s_3))$ where $F_X =\{ [s_1], [s_2], [s_3]\}$.
By theorem 2, $m_Y$ is determined by
$\pi_{F_Y}(m_Y)$ $= (t_{m_Y}(s_1^+), t_{m_Y}(s_1^-), t_{m_Y}(s_4),
t_{m_Y}(s_5), t_{m_Y}(s_6), t_{m_Y}(s_7))$.
Finally, formula (1)
shows that $t_{m_X}(s_7) = t_{m_X}(s_1) t_{m_X}(s_2)
t_{m_X}(s_3) + 2 - t_{m_X}^2(s_1) - t_{m_X}^2(s_2) - t_{m_X}^2(s_3).$
Combining these and lemma 12, we obtain the following lemma.

\midspace{0.1cm}
\centerline{\epsfbox{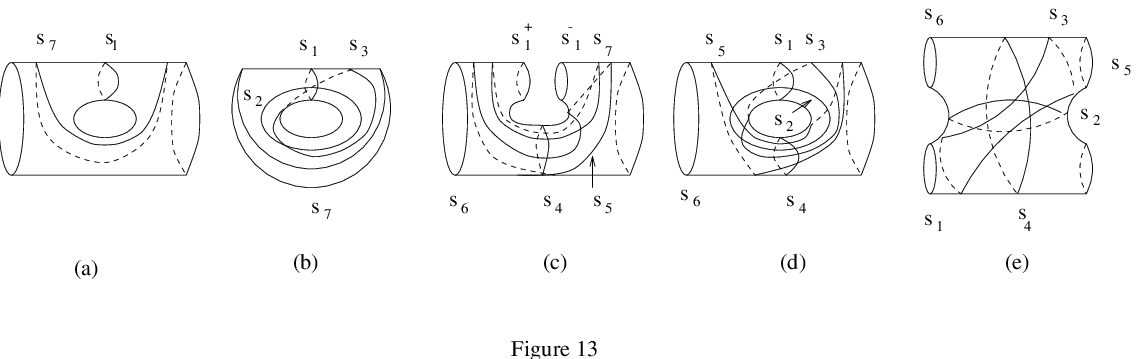}}
\midspace{0.1cm}

{\bf Lemma 13.} \it For surface $\Sigma_{1,2}$, let $F$ be the collection
of isotopy classes of six curves $s_1$, $s_2$, $s_3$, $s_4$, $s_5$, $s_6$
as in figure 13(d). Then
$\pi_F : T(\Sigma_{1,2}) \to R^6$ is an embedding whose image is given
by $\{(t_1, t_2, t_3, t_4, t_5, t_6) \in \bold R_{>2}^6 |$
$t_1t_2t_3 > t_1^2 + t_2^2 + t_3^3$ and 
$t_4 t_5 t_7 > t_4^2 + t_5^2 + t_6^2 + t_7^2 + 2t_1^2 + 2t_1^2t_6
+ t_1^2t_7 + t_1t_4t_6 + t_1t_5t_6 + 2t_1t_4 + 2t_1t_5 + 2t_6t_7$,
where $t_7 = t_1t_2t_3 -t_1^2 -t_2^2 -t_3^2$\}. Furthermore, there
exists an RC-section for $T(\Sigma_{1,2})$ with associated set $F$.
\rm

5.3. \it Proof of the Corollary for $\Sigma_{g,r}$ with $r >0$\rm

We prove the corollary  by induction on
$|\Sigma_{g,r}| = 3g+r$ with $r>0$.

For surfaces $\Sigma_{0,3}$, $\Sigma_{0,4}$, $\Sigma_{1,1}$ and
$\Sigma_{1,2}$, we have shown in the previous sections that
the corollary holds. 
Given $\Sigma_{g,r}$ with either $3g+r =n >5$ or $(g,r) = (0,5)$, 
if $r \geq 2$, we decompose $\Sigma_{g,r} = X \cup Y$
where $X \cong \Sigma_{g, r-1}$, $Y \cong \Sigma_{0,4}$ with $X \cap Y
 \cong \Sigma_{0,3}$ as in figure 3(b); if $r =1$, we decompose $\Sigma_{g,r}
=X \cup Y$ where $X \cong \Sigma_{g-1, 2}$, $Y \cong \Sigma_{1,2}$ and
$X \cap Y  \cong \Sigma_{0,3}$ as in figure 3(c). Then $|X|$ and $|Y|$ are
less than $|\Sigma_{g,r}|$. 
By the induction hypothesis, there exists a subset
$F_X \subset \Cal S(X)$ consisting of $6g+3r-9$ elements so that
corollary holds. Let $F_Y \subset \Cal S(Y)$ be the set
$\{ [s_1], [s_2],[s_3], [s_4], [s_5], [s_6]\}$
 given by theorem 2 as in figure 13(e) if
$Y \cong \Sigma_{0,4}$ and by lemma 13 as in figure 13(d) if 
$Y \cong \Sigma_{1,2}$.
Let $F = F_X \cup \{ [s_2], [s_3], [s_5]\}$ consisting of $6g+3r-6$ elements.
We claim that the corollary holds for $\Sigma_{g,r}$ 
with respect to the set $F$.
First to show that $\pi_F$ is an embedding, we use the gluing lemma. It
follows that $\pi_{F_X \cup F_Y}$ is an embedding. However, by the
construction,  $s_1, s_4 $ and $s_6$ are in the subsurface $X$. Thus
by the  induction hypothesis, $t_m(s_1)$, $t_m(s_4)$
and $t_m(s_6)$ are real analytic RC-functions in the coordinates of $\pi_{F_X}(m)$.
Thus, we may drop the three elements $[s_1], [s_4], $ and $[s_6]$ from
the set $F_X \cup F_Y$ without effecting the embeddedness of $\pi_{F_X
\cup F_Y}$. Applying lemma 12 to $F_X$ and $F_Y$ and then dropping
the three elements $[s_1], [s_4]$ and $[s_6]$, we see that $T(\Sigma_{g,r})$
has an RC-section with associated set $F$. Finally, we show that the
image $\pi_F(T(\Sigma_{g,r}))$ is defined by a finite set of
RC-inequalities in the coordinates of $\pi_F$.
Indeed, by the induction hypothesis, $\pi_{F_X}(T(X))$  (resp.
 $\pi_{F_Y}(T(Y))$) is
defined by a finite set of RC-inequalities.
By the gluing lemma 1, the image $\pi_{F_X \cup F_Y}(T(\Sigma_{g,r})$)
is given by the same set of RC-inequalities for $\pi_{F_X}(T(X))$,
together with the RC-inequalities for $\pi_{F_Y}(T(Y))$, and
three equations expressing that the lengths of the three simple closed
curves in $\partial ( X \cap Y)$ are the same in both metrics $m_X$ and
$m_Y$. Thus the result follows.

5.4. \it Proof of the Corollary  for closed surface $\Sigma_{g,0}$ with
$g \geq 2$\rm

Given $ \Sigma_g = \Sigma_{g,0}$, let $Y$ be an incompressible subsurface
of $\Sigma_g$ homeomorphic to $\Sigma_{1,1}$ with boundary $s_1$
and let $s_2$ be a non-separating simple closed curve in $int(Y)$. Set
$X = \Sigma_g -s_2$ as in figure 3(d). Thus $\Sigma_g = X \cup Y$ and $X \cap Y = Y -s_2$.
By the gluing lemma 1, each metric $m \in T(\Sigma_g)$ is the same as a
pair $(m_X, m_Y) \in T(X) \times T(Y)$ with $R_{X \cap Y}(m_X) = R_{X \cap Y}
(m_Y)$. In particular the the completion $\bar X$ of $X$ under the metric
$m_X$ has the same geodesic lengths at the two boundary components.
The following lemma describes hyperbolic
metrics on $\Sigma_{0,4}$ which have the same 
lengths at two boundary curves.

{\bf Lemma 14.} \it Given $\Sigma_{0,4}$ with curves $b_i$ ($i =1,2,3,4$)
as boundary components, let $a_{ij}$ $((i,j) = (1,2), (2,3),$ $(3,1))$
be  simple closed curves in $\Sigma_{0,4}$ so that $a_{12} \perp_0
a_{23}$ and $a_{31} = a_{12} a_{23}$ and $b_i$, $b_j$ and $a_{ij}$ bound
a subsurface of signature $(0,3)$.
Let $T'(\Sigma_{0,4})$ be the subspace of the
Teichm\"uller space $T(\Sigma_{0,4})$ defined by $t_m(b_3) = t_m(b_4)$,
and let $F' = \{[b_1], [b_2], [a_{12}], [a_{23}], [a_{31}]\}$. Then
$\pi_{F'} : T'(\Sigma_{0,4}) \to  \bold R_{>2}^5$ is an embedding whose image
is defined by a  real analytic RC-inequality in the coordinates of $\pi_{F'}$.
Furthermore, there is an RC-section $f$$: T'(\Sigma_g) \to \tilde R'(\Sigma_g)$
where $\tilde R'(\Sigma_g)$ stands for the subset of $\tilde R(\Sigma_g)$
which projects onto $T'(\Sigma_g)$ so that the entries of $f(m)(\gamma)$
are real analytic RC-functions in the coordinates of $\pi_{F'}(m)$. \rm

{\bf Proof}. Given a metric $m \in T'(\Sigma_g)$, let $t_i = t_m([b_i])$,
$i=1,2,3,4$, and let $t_{ij} = t_m([a_{ij}])$, $(i,j) = (1,2), (2,3)$, $(3,1)$,
where $t_3 = t_4$. Now these $t_i$ and $t_{ij}$ satisfy the equation
(13). Thus we obtain an equation in $t$ ($= t_3 = t_4$) below,
$$(2 + t_1t_2 + t_{12}) t^2 + (t_1t_{31} +  t_1t_{23} + t_2t_{31} + t_2t_{23}) t
+ t_1^2 + t_2^2 + t_1t_2t_{12}+ t_{12}^2 + t_{23}^2 + t_{31}^2 - t_{12}t_{23}t_{31} -4 = 0.$$

The coefficient of $t^2$ is positive and the constant term is negative
by (4). Thus the equation has two real roots of different signs and
$t_3$ ($ = t_4$) is the positive root of the equation. Thus $t_3 (= t_4$) is a
real analytic 
RC-function of $t_1$, $t_2$, $t_{12}$, $t_{23}$ and $t_{31}$ which
are the coordinates of $\pi_{F'}(m)$. This shows that $\pi_{F'}$ is
an embedding. The rest of the lemma follows by the same argument used
in the proof of theorem 2.
$\square$

Let $T'(X)$ be the subset of $T(X)$ so that $t_m(s_2^+) = t_m(s_2^-)$
where $s_2^+$ and $s_2^-$ are the boundary components of $\bar X$. Then
in the proof of the corollary for $\Sigma_{g-1, 2}$ ($\cong X$) in
\S5.3, to construct $m \in T'(X)$,
we decompose $X = X_1 \cup Y_1$ where $X_1 \cong \Sigma_{g-1, 1}$,
$Y_1 \cong \Sigma_{0,4}$ and $X_1 \cap Y_1 \cong \Sigma_{0,3}$.
We use lemma 14 instead of theorem 2 for metrics on $Y_1$ in the
gluing process. Thus, the same argument shows that
there exists a subset $F_X \subset \Cal S(\bar X)$
consisting of $6g-7$ elements so that $\pi_{F_X} : T'(X) \to
\bold R_{>2}^{6g-7}$ is an embedding whose image is an open 
set defined by a finite set of  real analytic RC-inequalities
in the coordinates of $\pi_{F_X}$. 

Let $s_3$ and $s_4$ be two simple closed curves in $int(Y)$ so that
$s_3 \perp s_2$ and $s_4 = s_2 s_3$.
Now by the gluing lemma 1, each $m \in$
 $T( X \cup Y)$ is determined by a pair $(m_X, m_Y) \in T'(X) \times T(Y)$
so that the restrictions of $m_X$ and $m_Y$ to $X \cap Y$ are the same.
The gluing condition on $X \cap Y$ is equivalent to that $t_{m_X}(s_2^+)$ 
$= t_{m_Y}(s_2)$ and $t_{m_X}(s_1) = t_{m_Y}(s_1) $ by lemma 4.
Also lemma 5 gives the complete description of $(t_{m_Y}(s_2), t_{m_Y}(s_3),$$
t_{m_Y}(s_4))$. Let $F = F_X \cup \{ [s_3], [s_4]\}$
$\subset \Cal S(\Sigma_g)$ consisting of $6g-5$ elements.
Combining the previous facts, we obtain (1) $\pi_F : T(\Sigma_g) \to
\bold R^{6g-5}$ is an embedding, (2) the image $\pi_F(T(\Sigma_g))$
is defined by a finite set of RC-inequalities (from those of $\pi_{F_X}
(T'(X))$ and of $\pi_{\{[s_2], [s_3], [s_4]\}}(T(\Sigma_{1,1}))$ 
where we replace $t_{m_Y}([s_2])$ by $t_{m_X}([s_2])$, and one 
real analytic RC-equation
$t_{m_Y}(s_1) = t_{m_X}(s_1)$.  
Furthermore, by lemmas 13 and  14, there is an RC-section for $T(\Sigma_g)$.
$\square$

\it Remark\rm.   The fact that $\pi_{F'}$ is an embedding in the lemma
10 was first proved by P. Schmutz ([Sc]).

\centerline{\bf Reference}

[Be] Bers, L.: Spaces of Riemann surfaces. Proc. Int. Congr. Math. Cambridge
(1959), 349-361

[Bi] Birman, J.: Mapping class groups of surfaces. In: Birman, J.,
Libgober, A. (eds.), Braids. Proceedings of a summer research
conference, Contemporary Math. Vol. 78, pp.13-44, Amer. Math. Soc. 1988

[Bo1] Bonahon, F.: The geometry  of Teichm\"uller space via geodesic
currents. Invent. Math. {\bf 92} (1988), 139-162

[Br] Brumfiel, G. W.: The real spectrum compactification of
Teichm\"uller space, Contemp. Math., {\bf 74}, AMS, (1988), 51-75

[Bu] Buser, P.: Geometry and spectra of compact Riemann surfaces. Progress
in Mathematics. Birkh\"auser, Boston, 1992

[CB] Casson, A., Bleiler, S.: Automorphisms of surfaces after Nielsen and
Thurston. Lond. Math. Soc. Student Text {\bf 9}. Cambridge University Press,
1988

[CS] Culler, M., Shalen, P.: Varieties of group representations and
splittings of 3-manifolds. Ann. Math. {\bf 117} (1983), 109-146

[De] Dehn, M.: Papers on group theory and topology. J. Stillwell (eds.).
 Springer-Verlag, Berlin-New York, 1987

[FK] Fricke, R., Klein, F.: Vorlesungen  \"uber die Theorie der
Automorphen Functionen. Teubner, Leipizig, 1897-1912

[FLP] Fathi, A., Laudenbach, F., Poenaru, V.: Travaux de Thurston sur les
surfaces. Ast\'erisque {\bf 66-67}, Soci\'et\'e Math\'ematique de France, 1979

[GM] Gilman, J., Maskit, B.: An algorithm for two-generator discrete groups.
Mich. Math. J. {\bf 38} (1) (1991), 13-32

[Go1] Goldman, W.: Topological components of spaces of representations.
Invent. Math. {\bf 93}(3) (1988), 557-607

[Go2] Goldman, W.: Introduction to character varieties. Preprint.

[Har] Harer, J.: The second homology group of the mapping class group
of an orientable surface. Invent. Math. {\bf 72} (1983), 221-239

[Hav] Harvey, W.: Spaces of discrete groups. In: Harvey, W. (ed.) Discrete
groups and automorphic functions. Proceedings of a NATO Advanced Study
Institute, pp.295-348. Academic Press, New York, 1977

[He] Helling, Diskrete Untergruppen von SL(2, $\bold R$), Inventiones Math.,
{\bf 17}, (1972), 217-229

[Ho] Horowitz, R.: Characters of free groups represented in the 2-dimensional
special linear group. Comm. Pure and App. Math. {\bf 25} (1972), 635-649

[HT] Hatcher, A. Thurston, W.: A presentation for the mapping class group
of a closed orientable surface. Topology {\bf 19} (1980), 221-237

[Jo] Johnson, D.: Spin structures and quadratic forms on surfaces. J. Lond.
Math. Soc. {\bf 22} (2) (1980), 365-373

[Jor] J$\phi$rgensen, T.: Closed geodesics on Riemann surfaces. 
Proc. Amer. Math.  Soc. {\bf 72} (1978), 140-142

[Ke] Keen, L.: Intrinsic moduli on Riemann surfaces. Ann. Math. {\bf 84} (1966),
405-420

[Li] Lickorish, R.: A representation of oriented combinatorial 3-manifolds. Ann.
Math. {\bf 72} (1962), 531-540

[Lu1] Luo, F.: Non-separating simple closed curves in a compact 
surface. Topology, {\bf 36} (2),  (1997), 381-410.

[Lu2] Luo, F.: Simple loops on surfaces and their intersection numbers, preprint.

[Lu3] Luo, F.:  Geodesic length functions and Teichm\"uller spaces,
Electronic Research Announcement, AMS, {\bf 2} (1), (1996), 34-41

[Mag] Magnus, W.: Rings of Fricke characters and automorphism groups of
free groups. Math. Zeit. {\bf 170} (1980), 91-103

[Mas] Maskit, B.: On Klein's combination theorem. Trans. Amer. Math. Soc.
{\bf 131} (1968), 32-39

[MSh] Morgan, J., and Shalen, P.: Valuations, trees, and degenerations of
hyperbolic structures, I, Ann. of Math., {\bf 120} (1984), 401-476

[MSi] Macbeth, A. M., Singerman, D.: Spaces of subgroups and 
Teichm\"uller space. Proc. Lond. Math. Soc. {\bf 31} (3) (1975), 211-256

[Ok1] Okumura, Y.: On the global real analytic coordinates for Teichm\"uller
spaces. J. Math. Soc. Jap. {\bf 42} (1990), 91-101

[Ok2] Okumura, Y.: Global real analytic length parameters for
Teichm\"uller spaces. Hiroshima Math. J.  {\bf 26} (1) (1996), 165--179

[Sc] Schmutz, P.: Die Parametrisierung des Teichm\"ullerraumes durch
geod\"atishe L\"angenfun \newline ktionen. Comm. Math. Hel. {\bf 68} (1993), 278-288

[SS] Sepp\"al\"a, M and Sorvali, T.: On geometric parametrization of 
Teichm\"uller spaces, Ann. Acad. Sci. Fenn. Ser. A I Math., {\bf 10}
(1985), 515-526

[So] Sorvali, T.: Parametrization for free M\"obius groups. Ann. Acad. Sci.
Fenn. {\bf 579} (1974), 1-12

[Th] Thurston, W.: On the geometry and dynamics of diffeomorphisms of
surfaces. Bul. Amer. Math. Soc. {\bf 19} (2) (1988), 417-438       

[Wo] Wolpert, S.: Geodesic length functions and the Nielsen problem.
J. Diff. Geo. {\bf 25} (1987), 275-296

% Department of Mathematics Rutgers University New Brunswick, NJ 08903 

\end